\documentclass[11pt]{amsart}

\title[Hopf Algebras in Model Categories to Hopf Algebras in \texorpdfstring{$\infty$}{∞}-Categories]{From Hopf Algebras in Model Categories to Hopf Algebras in \texorpdfstring{$\infty$}{∞}-Categories}

\author[Klang]{Inbar Klang}
\address[Klang]{
Department of Mathematics,
Vrije Universiteit Amsterdam - Faculty of Science,
De Boelelaan 1111,
1081 HV Amsterdam,
The Netherlands}
\email{i.klang@vu.nl}

\author[Kuijper]{Josefien Kuijper}
\address[Kuijper]{
Department of Mathematics, University of Toronto
Bahen Centre, Room 6290
40 St. George St.,
Toronto, Ontario M5S 2E4,
Canada}
\email{josefien.kuijper@utoronto.ca}

\author[Malkiewich]{Cary Malkiewich}
\address[Malkiewich]{
Department of Mathematics, Binghamton University,
PO Box 6000, 
Binghamton, New York 13902,
USA}
\email{cmalkiew@binghamton.edu}

\author[Mehrle]{David Mehrle}
\address[Mehrle]{
{Department of Mathematics, 
Colorado State University, 
Louis R. Weber Building, 
841 Oval Drive, 
Fort Collins, Colorado 80523
USA}}
\email{{david.mehrle@colostate.edu}}

\author[Wittich]{Thor Wittich}
\address[Wittich]{
Department of Mathematics, Universit\"{a}t Osnabr\"{u}ck, 
Albrechtstrasse 28a, 49076 Osnabr\"{u}ck,
Germany}
\email{thor.wittich@uni-osnabrueck.de}

\usepackage{amssymb,amsthm,amsmath} 
\usepackage{bbm} 
\usepackage{enumitem} 
\usepackage[margin=1.25in]{geometry} 
\usepackage[pagebackref,
			colorlinks, 
			citecolor=OI3, 
			linkcolor=OI3, 
			urlcolor=OI5]{hyperref} 
\usepackage{mathrsfs} 
\usepackage{thmtools} 
\usepackage{tikz} 
\usepackage{tikz-cd} 
\usepackage{todonotes} 
\usepackage{xcolor} 
\usepackage{adjustbox}
			
\usepackage[nameinlink,capitalize]{cleveref} 

\definecolor{OI1}{RGB}{230,159,0}	
\definecolor{OI2}{RGB}{86,180,233}	
\definecolor{OI3}{RGB}{0,158,115}	
\definecolor{OI4}{RGB}{240,228,66}	
\definecolor{OI5}{RGB}{0,114,178}	
\definecolor{OI6}{RGB}{213,94,0}	
\definecolor{OI7}{RGB}{204,121,167}	

\DeclareMathAlphabet\mathbfcal{OMS}{cmsy}{b}{n} 

\DeclareMathOperator{\colim}{colim} 
\DeclareMathOperator{\Ex}{Ex} 
\newcommand{\id}{\mathrm{id}} 
\DeclareMathOperator{\sh}{sh} 

\DeclareMathOperator{\op}{\textnormal{op}}
\DeclareMathOperator{\lax}{\textnormal{lax}}

\DeclareMathOperator{\Alg}{Alg} 
\DeclareMathOperator{\Assoc}{Assoc} 
\DeclareMathOperator{\BiAlg}{BiAlg} 
\DeclareMathOperator{\CAlg}{CAlg} 
\DeclareMathOperator{\Cat}{Cat} 
\DeclareMathOperator{\CBiAlg}{CBiAlg} 
\DeclareMathOperator{\CoAlg}{CoAlg} 
\DeclareMathOperator{\Comm}{Comm} 
\DeclareMathOperator{\CHopf}{CHopf} 
\DeclareMathOperator{\Fin}{Fin} 
\newcommand{\Fun}{\mathrm{Fun}} 
\DeclareMathOperator{\sCat}{sCat} 
\DeclareMathOperator{\Sp}{Sp} 
\DeclareMathOperator{\Spp}{\mathbf{Sp}} 
\DeclareMathOperator{\sSet}{sSet} 

\makeatletter
\newcommand*{\boxwedge}{\mathbin{\mathpalette\@boxwedge{}}}
\newcommand*{\@boxwedge}[2]{%
  \sbox0{$#1\boxtimes\m@th$}%
  \dimen2=.5\dimexpr\wd0-\ht0-\dp0\relax 
  \dimen@=\dimexpr\ht0+\dp0\relax
  \def\lw{.07} 
  \kern\dimen2 
  \tikz[
    line width=\lw\dimen@,
    line join=round,
    x=\dimen@,
    y=\dimen@,
    baseline=0
  ]
  \draw
    (\lw/2,0) rectangle (1-\lw,1-\lw)
    (\lw,0) -- (.5,1-\lw-\lw/2) -- (1-\lw-\lw/2,0)
  ;%
  \kern\dimen2 
}
\makeatother

\numberwithin{equation}{section}
\numberwithin{figure}{section}

\newtheorem{corollary}[equation]{Corollary}
\AddToHook{env/corollary/begin}{\crefalias{equation}{corollary}}
\newtheorem{lemma}[equation]{Lemma}
\AddToHook{env/lemma/begin}{\crefalias{equation}{lemma}}
\newtheorem{proposition}[equation]{Proposition}
\AddToHook{env/proposition/begin}{\crefalias{equation}{proposition}}
\newtheorem{theorem}[equation]{Theorem}
\AddToHook{env/theorem/begin}{\crefalias{equation}{theorem}}

\newtheorem{lettertheorem}{Theorem}
\AddToHook{env/lettertheorem/begin}{\crefalias{lettertheorem}{theorem}}
\newtheorem{lettercorollary}[lettertheorem]{Corollary}
\AddToHook{env/lettercorollary/begin}{\crefalias{lettertheorem}{corollary}}

\theoremstyle{definition}
\newtheorem{construction}[equation]{Construction}
\AddToHook{env/construction/begin}{\crefalias{equation}{construction}}
\newtheorem{definition}[equation]{Definition}
\AddToHook{env/definition/begin}{\crefalias{equation}{definition}}
\newtheorem{example}[equation]{Example}
\AddToHook{env/example/begin}{\crefalias{equation}{example}}
\newtheorem{remark}[equation]{Remark}
\AddToHook{env/remark/begin}{\crefalias{equation}{remark}}
\newtheorem{warning}[equation]{Warning}
\AddToHook{env/warning/begin}{\crefalias{equation}{warning}}

\crefname{equation}{}{} 

\begin{document}


\begin{abstract}
	We show that algebra objects in model categories can be transferred to algebra objects in $\infty$-categories, without any cofibrancy or fibrancy assumptions on the algebra. 
	We furthermore show under some mild extra assumptions that this correspondence extends to commutative bialgebras and to commutative Hopf algebras.
\end{abstract}

\maketitle

\begingroup%
	\setlength{\parskip}{0em} 
	\setcounter{tocdepth}{1}
	\tableofcontents
\endgroup%

\section{Introduction}

The theory of $\infty$-categories has become an invaluable tool in mathematics, particularly in homotopy theory.
One reason for its incredible utility is its ability to encode homotopy coherence, which enables us to endow objects, categories, and maps with structure in a more flexible way than with earlier point-set approaches. 
For example, strict coalgebras in a symmetric monoidal model category of spectra are always cocommutative \cite{PS_coalgebras}, whereas coalgebras in the $\infty$-category of spectra are more plentiful and capture more examples of interest.

In order to make full use of the power of $\infty$-categories, it is essential to have frameworks for translating results from more concrete approaches (e.g., model categories) into $\infty$-categories. 
Many translational results have thus far been established. For instance, there is an equivalence between cofibrant algebras over a $\Sigma$-cofibrant operad $\mathcal O$ in a symmetric monoidal model category, and algebras over the corresponding $\infty$-operad in the underlying symmetric monoidal $\infty$-category \cite{HA,PavlovScholbach}, \cite[Theorem 4.10]{haugseng}, \cite[Theorem 7.3.3]{whiteyau}. 

One situation where these existing translational results do not suffice occurs in \cite{KKMMW1}, in which the authors endow reduced spherical scissors congruence $K$-theory with the structure of a Hopf algebra spectrum, using model categories. To get the same structure in the $\infty$-category of spectra, one needs a translation as above, but the algebra in question is not cofibrant, and any cofibrant replacement as an algebra would destroy its coalgebra structure.

In this paper, we remedy this by proving a translational result that does not make cofibrancy assumptions. 

\begin{lettertheorem}[\cref{body_theorem_A}]
\label{thm-alg-interchange-intro}
   Let $M_\bullet$ be a symmetric monoidal simplicial model category with underlying symmetric monoidal $\infty$-category $\mathbfcal{M}$, and let $\mathcal{O}_\bullet$ be a fibrant simplicial operad. 
   Then there is a canonical ``interchange'' map of simplicial sets
	\[ 
		N^s(\Alg_{\mathcal{O}_\bullet}(M_\bullet)) \to \Alg_{N^s(\mathcal{O}_\bullet)}(\mathbfcal{M}). 
	\]
\end{lettertheorem}

As a result, any algebra or diagram of algebras in $M_\bullet$ becomes an algebra or diagram of algebras in the $\infty$-category $\mathbfcal{M}$. 
This is true even if the operad $\mathcal O$ and the algebras over $\mathcal O$ are not cofibrant.

Note that in the statement of \cref{thm-alg-interchange-intro}, the simplicial set on the right is always an $\infty$-category. 
If we assume that $M_\bullet$ comes from a topological model category, then the simplicial set on the left is an $\infty$-category as well (\cref{top_case_1}), making this a functor of $\infty$-categories. 
However, even without this assumption, the simplicial set on the left describes a large and interesting class of ``concrete'' $\mathcal{O}_\bullet$-algebras in the underlying $\infty$-category $\mathbfcal{M}$.

\cref{thm-alg-interchange-intro} fills a marked gap in the existing literature, because the earlier work in \cite{HA,PavlovScholbach,haugseng,whiteyau} makes stronger assumptions in order to ensure that the above is not just a functor but an equivalence of $\infty$-categories. 
In the example of interest to us, this functor is not an equivalence, and we only need the one-way translation: every algebra in the model category creates an algebra in the underlying $\infty$-category. 

We also prove a corresponding result for commutative bialgebra and commutative Hopf algebra objects. 
In this case we make the technical assumption that $M_\bullet$ arises as the singular simplices of a symmetric monoidal topological model category. 
We also impose a cofibrancy condition, but a mild one. More precisely, our second main result is the following.

\begin{lettertheorem}[\cref{main_comparison}]\label{main_comparison-intro}
	Let $M_{\bullet}$ be the singular simplices of a symmetric monoidal topological model category, and let $M_{\bullet}^\flat$ be any full subcategory of $M_{\bullet}$ containing the cofibrant objects, closed under tensor product, and on which the tensor product preserves all equivalences. 
	Let $D_1$ be the little intervals operad, regarded as an operad in simplicial sets by taking singular simplices. 
	Then each $(\Comm,D_1)$-bialgebra in the subcategory $M^\flat_\bullet$ induces a commutative bialgebra in the underlying symmetric monoidal $\infty$-category $\mathbfcal{M}$ of $M_\bullet$. 
	This extends to a map of $\infty$-categories
	\[ 
		N^s(\BiAlg_{(\Comm,D_1)}(M^\flat_\bullet)) \to  \CBiAlg(\mathbfcal{M}) 
	\]
	which preserves the product, coproduct, and shear map as maps in the homotopy category $hM^\flat_\bullet$. 
	In particular, this map preserves Hopf algebra objects.
\end{lettertheorem}

The reason for the cofibrancy assumption in \cref{main_comparison-intro} is that without it, the passage from $M_\bullet$ to $\mathbfcal{M}$ would only be lax monoidal, and not strong monoidal (see \cref{underlying_sm_inf_cat_compare_2}). 
This would be enough to preserve the algebra structure, as in \cref{thm-alg-interchange-intro}, but not enough to preserve the coalgebra structure as well.

A special case of \cref{main_comparison-intro}, which we use in \cite{KKMMW1}, is the following. Let $I$ be a symmetric monoidal 1-category, and let $\Sp^O_\bullet$ denote the symmetric monoidal simplicial model category of orthogonal spectra. Let $\Fun(I,\Sp^O_\bullet)^\flat$ denote the 1-category of pointwise cofibrant diagrams of orthogonal spectra. Suppose $X:I\to \Sp^O_\bullet$ is a pointwise cofibrant commutative bialgebra for the Day convolution in $\Fun(I,\Sp^O_\bullet)^\flat$. Then with \cref{main_comparison-intro}, we can transfer $X$ to a commutative bialgebra in the symmetric monoidal $\infty$-category $\Fun(N(I),\mathbf{Sp})$ of diagrams in the $\infty$-category $\mathbf{Sp}$ of spectra. By applying the homotopy colimit, we obtain a commutative bialgebra in the $\infty$-category $\mathbf{Sp}$. 

\begin{lettercorollary}[\cref{main_corollary}]\label{main_corollary-intro}
	For any symmetric monoidal indexing $1$-category $I$, and for each $(\Comm,D_1)$-bialgebra $X$ in the symmetric monoidal model category $\Fun(I,\Sp^O_\bullet)$ that is pointwise cofibrant, the homotopy colimit $\mathrm{hocolim}\, X$ is a commutative bialgebra in the $\infty$-category of spectra $\mathbf{Sp}$. This extends to a map of $\infty$-categories
	\[ 
		N^s(\BiAlg_{(\Comm,D_1)}(\Fun(I,\Sp^O_\bullet)^\flat)) \to  \CBiAlg(\mathbf{Sp})
	\]
	which preserves the product, coproduct, and shear map as maps in the homotopy category, and therefore also preserves Hopf algebra objects.
\end{lettercorollary}

\begin{remark}\label{rem-comm-D1}
        \cref{main_comparison-intro} {and \cref{main_corollary-intro} pertain} specifically to $(\Comm,D_1)$-bialgebras, as this is the situation that arises in \cite{KKMMW1}. 
        The use of the commutative operad is significant in our argument, as the symmetric monoidal product in the category of commutative algebras is also the coproduct in this category. 
        However, the choice of the little intervals operad $D_1$ is not essential and corresponding result{s hold} for $(\Comm,\mathcal O_\bullet)$-bialgebras for any other fibrant simplicial operad $\mathcal{O}_\bullet$
        (except that the statements about Hopf algebras only make sense if coalgebras for $\mathcal{O}_\bullet$ create coassociative coalgebras in the homotopy category).
\end{remark}

\subsection*{Notation and Conventions}

We make the following conventions. 
We use ``$\infty$-categories" in this paper to refer to the quasicategories of Joyal and Lurie. 
We try to be explicit whether we are talking about $1$-categories, simplicial categories, or $\infty$-categories. 
The unqualified term ``category" shouldn't appear here. 
\begin{center}
	\begin{tabular}{ l l l }
		$1$-categories and constant simplicial categories & & $C, D, \dotsc$\\
		Simplicial categories & & $C_\bullet, D_\bullet, \dotsc$\\
		$\infty$-categories & & $\mathbfcal{C},\mathbfcal{D},\dotsc$\\
		Nerves of $1$-categories & & $N(C),N(D),\dotsc$\\
		Simplicial nerves (of simplicial categories) & & $N^s(C_\bullet),N^s(D_\bullet),\dotsc$\\
		Simplicial operads  & & $\mathcal{O}_\bullet,\mathcal{O}_\bullet',\dotsc$\\
		Symmetric monoidal $\infty$-categories  & & $\mathbfcal{C}^\otimes,\mathbfcal{D}^\otimes,\dotsc$\\
		$\infty$-operads/$\infty$-multicategories  & & $\mathbfcal{O}^\otimes,\mathbfcal{O}'^\otimes,\dotsc$\\
		Simplicial category of algebras over $\mathcal{O}_\bullet$ in $C_{\bullet}$ & & $\Alg_{\mathcal{O}_\bullet}(C_\bullet)$\\
		$\infty$-category of algebras over $\mathbfcal{O}^\otimes$ in $\mathbfcal{C}^\otimes$ & & $\Alg_{\mathbfcal{O}}(\mathbfcal{C})$\\
	\end{tabular}
\end{center}

\subsection*{Acknowledgments.}

The authors would like to thank the organizers of the Collaborative Research Workshop on K-theory and Scissors Congruence at Vanderbilt University in July 2024, during which this paper was initiated, and the conference Scissors congruence and K-theory at the University of Pennsylvania in July 2025, during which this paper was almost finished. 
Both were sponsored by the NSF grant FRG: Collaborative Research: Trace Methods and Applications for Cut-and-Paste K-Theory. 
The authors are grateful to Maximilien P\'{e}roux for extensive guidance on model categories of coalgebras, and how to ultimately avoid using them.
Additional thanks are due to Ulrich Bunke and Tobias Lenz for invaluable conversations on the technical points of localizations, and to Marco Giustetto and Maxime Ramzi for fruitful discussions about $\infty$-categorical technicalities.

JK was partially supported by the Knut and Alice Wallenberg grant KAW 2023.0416. 
CM was partially supported by the National Science Foundation (NSF) grants DMS-2052923 and DMS-2506430 and by a Simons Fellowship.
DM was partially supported by NSF grant DMS-2135884.
TW would like to thank the research training group 2240: Algebro-geometric Methods in Algebra, Arithmetic and Topology for partial support. 

Finally, we would like to thank the anonymous referee for their careful reading of this paper, and for their many helpful comments which improved the presentation of this paper.

\section{Symmetric monoidal \texorpdfstring{$\infty$}{∞}-categories}

We begin with a review of $\infty$-categories and symmetric monoidal $\infty$-categories, following \cite{HTT,HA}. 
We provide some details in those places where our later proofs require it.

\subsection{Definitions}

\begin{definition}\label{def:infty_category}
    An \textit{$\infty$-category} is a quasicategory, ie,  a simplicial set in which every inner horn has an extension \cite[1.1.2.4]{HTT}. 
    A \emph{functor} between $\infty$-categories is a morphism between the underlying simplicial sets. {A functor of $\infty$-categories is an \emph{equivalence} when it is a weak equivalence in the Joyal model structure from \cite[2.6]{joyal_notes} or \cite[2.2.5.1]{HTT}.}
\end{definition}

    For any 1-category $C$, let $N(C)$ denote its nerve, which is always an $\infty$-category. 
    For any simplicially enriched category $C_\bullet$, let $N^s(C_\bullet)$ denote its simplicial nerve or homotopy coherent nerve \cite[1.1.5.5]{HTT}.
    This is always a simplicial set, and it is an $\infty$-category when $C_\bullet$ is \emph{fibrant}, meaning that its mapping spaces $C_\bullet(a,b)$ are Kan complexes. {This follows from the fact that $N^s$ is a right Quillen functor \cite[Section 1.3]{Joyal2007} for the Bergner model structure on simplicial categories from \cite{bergner2007}; see also \cite[1.1.5.10]{HTT}.}
    
    We also have $N^s(C_\bullet) \cong N(C_0)$ when $C_\bullet$ is \emph{discrete}, meaning that the mapping spaces $C_\bullet(a,b)$ are all isomorphic to constant simplicial sets.

\begin{remark}
	\label{simplicial categories as sCat}
    We recall that a simplicially enriched category can equivalently be seen as a simplicial object in $1$-categories, where the $1$-categories $C_0, C_1,\dots$ have the same set of objects, and the face and degeneracy maps are functors which are constant on objects. 
    From now on, we will refer to these as \emph{simplicial categories}. 
\end{remark}

Let $(C,\otimes,\mathbbm{1})$ be a \emph{symmetric monoidal $1$-category}. 
This consists of the data of a $1$-category $C$ together with a tensor product functor $- \otimes - \colon C \times C \rightarrow C$, a unit object $\mathbbm{1} \in C$, and natural isomorphisms encoding the associativity, unitality, and commutativity of the tensor product, such that certain diagrams commute. 
The most na\" ive translation of this definition to $\infty$-categories results in a definition that is so cumbersome that it is unusable. 
Instead, we encode the symmetric monoidal structure using Grothendieck opfibrations, as we recall next.

For each integer $n \geq 0$, let $[n]_+$ be the finite pointed set $\{1,2,\ldots, n,*\}$ with basepoint $*$. 
More generally, for any finite set $I$, let $I_+$ be $I$ with a disjoint basepoint.

\begin{definition}\label{category_of_operators}
	For any symmetric monoidal 1-category $(C,\otimes,\mathbbm{1})$, we define $C^{\otimes}$ to be the 1-category whose objects are tuples of objects $\{A_i\}_{i \in I}$, where $I$ is allowed to be any finite set. 
	A morphism $\{A_i\}_{i \in I} \to \{B_j\}_{j \in J}$ in $C^\otimes$ is given by a pair $(\alpha,\{f_j\}_{j \in J})$ where $\alpha \colon I_+ \to J_+$ is a basepoint-preserving function and each 
	\[
		f_j \colon \bigotimes_{i \in \alpha^{-1}(j)} A_i \to B_j
	\]
	is a morphism in $C$. 
	Composition in $C^\otimes$ composes the functions and the morphisms in $C$:
	\[
		\big(\beta,\{g_k\}\big) \circ \big(\alpha, \{f_j\}\big) 
			= \bigg(\beta \circ \alpha, \Big\{g_k \circ \bigotimes_{j \in \beta^{-1}(k)} f_j \Big\}\bigg).
	\]
\end{definition}

If we let $\Fin_*$ denote the $1$-category of finite pointed sets, then there is a forgetful functor $p \colon C^{\otimes} \rightarrow \Fin_*$ that sends each tuple $\{A_i\}_{i \in I}$ to the set $I_+$, and for each morphism only remembers the function $\alpha\colon I_+ \to J_+$. 
We denote the fiber of $p$ over $I_+$ by $C^\otimes_{I_+}$.

Note that the functor $p$ is an \emph{opfibration}, meaning that for each basepoint-preserving function $\alpha\colon I_+ \to J_+$ and tuple $\{A_i\}_{i \in I}$, there is a privileged \emph{cocartesian morphism}
\[ 
	\bar\alpha\colon \{A_i\}_{i \in I} \to \Big\{ \bigotimes_{i \in \alpha^{-1}(j)} A_i \Big\}_{j \in J} 
\]
characterized by the property that composing with $\bar\alpha$ gives a bijection between those morphisms 
\(
	\{A_i\}_{i \in I} \to \{B_j\}_{j \in J}
\) 
that lie over $\alpha$, and morphisms 
\(
	\big\{ \bigotimes_{i \in \alpha^{-1}(j)} A_i \big\}_{j \in J} \to \{B_j\}_{j \in J}
\) 
that lie over the identity of $J$. 
More concretely, $\bar\alpha$ is the morphism that selects only identity maps of the tensors $\bigotimes_{i \in \alpha^{-1}(j)} A_i$. 
The cocartesian morphisms $\bar\alpha$ allow us to define \emph{pushforward} functors $\alpha_!\colon C^\otimes_{I_+} \to C^\otimes_{J_+}$, which are just the functors that take the tensor product along the preimages of $\alpha$.

Furthermore, $p$ satisfies the following Segal condition: if we take the $n$ different maps 
\(
	\rho^i\colon [n]_+ \to [1]_+
\) 
that fold all but one of the points into the basepoint, then the pushforwards $\rho^i_!$ collectively define an equivalence of $1$-categories 
\(
	C^\otimes_{I_+} \simeq \prod_{i \in I} C^\otimes_{{[}1{]}_+}.
\)

The structure on $C^\otimes$ that we have just discussed characterizes the symmetric monoidal structure on $C$, and it translates much better into the setting of $\infty$-categories:

\begin{definition}\label{df:sm_infinity_cat}
    A \emph{symmetric monoidal $\infty$-category} is an $\infty$-category $\mathbfcal{C}^{\otimes}$ and a functor $p \colon \mathbfcal{C}^{\otimes} \rightarrow N(\Fin_*)$ that is a cocartesian fibration (\cite[2.4.2.1]{HTT}, \cite[3.1]{bs_fibrations}) and that satisfies the Segal condition analogous to the one above \cite[Def 2.0.0.7]{HA}. 
    
    More generally, an \emph{$\infty$-multicategory}, also called an \textit{$\infty$-operad}, is an $\infty$-category $\mathbfcal{O}^{\otimes}$ and a functor $p \colon \mathbfcal{O}^{\otimes} \rightarrow N(\Fin_*)$ that satisfies a Segal condition and a condition that is weaker than being a cocartesian fibration \cite[Def 2.1.1.10]{HA}.
\end{definition}

In particular, a symmetric monoidal $\infty$-category has an underlying $\infty$-category $\mathbfcal{C}$, defined to be the fiber $\infty$-category $\mathbfcal{C}^\otimes_{[1]_+} = p^{-1}([1]_+)$, and tensor product and unit functors
\[ 
	\otimes\colon \mathbfcal{C} \times \mathbfcal{C} \to \mathbfcal{C}, \qquad I\colon * \to \mathbfcal{C}, 
\]
defined as the pushforward along the map $[2]_+ \to [1]_+$ that folds the two points into one, and the map $[0]_+ \to [1]_+$, respectively. 
Note that these functors are only well-defined up to natural isomorphism. 
The rest of the structure encodes associativity, commutativity, and unitality of the tensor product up to various isomorphisms, along with an infinite collection of coherences between these isomorphisms. 
In particular, we get a symmetric monoidal structure on the homotopy category $h\mathbfcal C$, whose product functor is well-defined up to canonical isomorphism.

The following consequence of \cite[2.4.1.10]{HTT} is useful for checking the condition of being a cocartesian fibration. 
For any map of simplicial categories $p\colon C_{\bullet} \to \Fin_*$, where we regard $\Fin_*$ as a constant simplicial category, let $(C_{\bullet})_{I_+}$ denote the fiber category over the object $I_+$, and for any objects $x$ and $y$ in $C_{\bullet}$, let $(C_{\bullet})_\alpha(x,y)$ denote the component of the mapping space $C_{\bullet}(x,y)$ lying over a given morphism $\alpha$ in $\Fin_*$.

\begin{lemma}\label{cocartesian_criterion}
    {Let $C_\bullet$ be a fibrant simplicial category.} A map of simplicial categories $p\colon C_{\bullet} \to \Fin_*$ becomes a cocartesian fibration of $\infty$-categories $N^s(C_{\bullet}) \to N(\Fin_*)$ if for every morphism $\alpha\colon I_+ \to J_+$ and every object $x \in (C_{\bullet})_{I_+}$, there is a morphism $\bar\alpha\colon x \to y$ over $\alpha$ that is ``cocartesian'' in the sense that for every morphism $\beta\colon J_+ \to K_+$, composing with $\bar\alpha$ gives a weak equivalence of Kan complexes
    \[ 
    	(- \circ \bar{\alpha})\colon (C_{\bullet})_{\beta}(y,z) \to (C_{\bullet})_{\beta\circ\alpha}(x,z). 
	\]
\end{lemma}

\subsection{Examples from symmetric monoidal simplicial categories}

Next we explain how to produce concrete examples of symmetric monoidal $\infty$-categories, using the simpler notion of a symmetric monoidal simplicial category.

\begin{definition}\label{sm_infinity_cat_most_elementary_example}
    A \emph{symmetric monoidal simplicial category} is a simplicial category $C_\bullet$ with a symmetric monoidal structure that acts on the simplicial enrichment. 
    	In other words, seen as a simplicial object in {1-}categories, each category $C_n$ is given a symmetric monoidal structure, such that the face and degeneracy maps are strict symmetric monoidal functors that are the identity on the set of objects.
    
    When $C_\bullet$ is a symmetric monoidal simplicial category, we form the simplicial category $C^\otimes_\bullet$ by applying \cref{category_of_operators} to each simplicial level. This comes with a morphism $C^\otimes_\bullet \to \Fin_*$ to the constant simplicial category $\Fin_*$.
\end{definition}

\begin{lemma}\label{ssm_to_smi}
    If $C_\bullet$ is a symmetric monoidal simplicial category and it is fibrant as a simplicial category, then $N^s(C^\otimes_\bullet) \to N(\Fin_*)$ is a symmetric monoidal $\infty$-category.
\end{lemma}

\begin{proof}
    The map $N^s(C^\otimes_\bullet) \to N(\Fin_*)$ is an $\infty$-operad by \cite[2.1.1.27]{HA}, so it has the Segal condition already. 
	By \cref{cocartesian_criterion}, it suffices to show that for every morphism $\alpha\colon I_+ \to J_+$ in $\Fin_*$ and object $\{A_i\}_{i \in I}$ in $C^\otimes$, there is a 0-cell $\bar{\alpha}$ in the space of morphisms, which lies over $\alpha$, such that composition with $\bar{\alpha}$ induces a weak equivalence between the maps lying over any given map $\beta\colon J_+ \to K_+$ and the maps lying over the composite $\beta \circ \alpha\colon I_+ \to K_+$.
	We choose $\bar{\alpha}$ just as we did after \cref{category_of_operators}. It induces an isomorphism on the mapping spaces, which is certainly a weak equivalence, and so $\bar{\alpha}$ is the desired cocartesian lift.
\end{proof}

Note that if $C_\bullet$ fails to be fibrant then we may fix this by applying Kan's fibrant replacement functor $\Ex^\infty$ to all of its mapping spaces. 
Since $\Ex^\infty$ preserves finite products up to isomorphism, we get:

\begin{lemma}\label{monoidal_fibrant_replacement}
	For any symmetric monoidal simplicial category $C_\bullet$, the simplicial category $\Ex^\infty C_\bullet$ is also a symmetric monoidal simplicial category, and $C_\bullet \to \Ex^\infty C_\bullet$ is a strong symmetric monoidal functor.
\end{lemma}
	
Therefore even though $N^s(C^\otimes_\bullet)$ is not a symmetric monoidal $\infty$-category, it maps forward to $N^s((\Ex^\infty C_\bullet)^\otimes)$, which is a symmetric monoidal $\infty$-category.

\begin{example}\label{cf_example}
    Suppose $M_\bullet$ is a symmetric monoidal simplicial category and it is also a model category, in which the tensor product preserves the cofibrant objects and the weak equivalences between them. Let $M_\bullet^{\mathrm{cf}}$ be the subcategory of objects that are \emph{bifibrant} (both cofibrant and fibrant). Then the construction from \cref{category_of_operators} for the subcategory of bifibrant objects gives a map of $\infty$-categories $N^s((M^{\mathrm{cf}}_\bullet)^\otimes) \to N(\Fin_*)$. The above proof does not show that this is a symmetric monoidal $\infty$-category, because the tensor product of bifibrant objects is only cofibrant, and not necessarily bifibrant. However, it is still a symmetric monoidal $\infty$-category by \cite[4.1.7.10]{HA}. We form the cocartesian arrows by taking the map to the tensor product and then composing with a fibrant replacement, giving a map whose target is indeed bifibrant.
\end{example}

\begin{definition}\label{underlying_sm_inf_cat_1}
    When $M_\bullet$ satisfies the assumptions of \cref{cf_example}, we say $N^s((M^{\mathrm{cf}}_\bullet)^\otimes)$ is the \emph{underlying symmetric monoidal $\infty$-category} of $M_\bullet$, following \cite[4.1.7.6]{HA}.
    
    We will see later in \cref{underlying_sm_inf_cat_2} that we get a better definition of ``underlying symmetric monoidal $\infty$-category'' by taking the entire $\infty$-category $N^s(M_\bullet)$ and then inverting all of the weak equivalences.
\end{definition}
    
\begin{example} 
	The underlying symmetric monoidal $\infty$-category of orthogonal spectra is $N^s(((\Sp^O)^{\mathrm{cf}}_\bullet)^\otimes)$. 
	This is one model for the symmetric monoidal $\infty$-category of spectra $\Spp$; see \cite[4.1.8.6]{HA}.   
\end{example}    

\begin{example}\label{ex:infty_op_from_simplicial_op}
	Let $\mathcal{O}_\bullet$ be any fibrant simplicial operad (meaning the simplicial sets $\mathcal O_\bullet(n)$ are Kan complexes). 
	We form the simplicial category $\mathcal O^\otimes_\bullet \to \Fin_*$ as in \cref{category_of_operators}.
	The objects of $\mathcal O_\bullet^\otimes$ are also the finite based sets, but the morphisms are the products of levels of the operad
    \[ 
    	\mathcal O_\bullet^\otimes(I_+,J_+) 
			:= 
		\coprod_{\alpha\colon I_+ \to J_+} \left( \prod_{j \in J} \mathcal O_\bullet(\alpha^{-1}(j)) \right), 
	\]
    and the compositions arise from the composition in $\mathcal O_{\bullet}$. 
    Since the simplicial sets $\mathcal O_\bullet(n)$ are Kan complexes, the simplicial nerve $N^s(\mathcal O_\bullet^\otimes) \to N(\Fin_*)$ is an $\infty$-operad by \cite[2.1.1.27]{HA}.
    
    As a special case, if $\mathcal O_\bullet$ is the singular simplices of the little intervals operad, which we call $D_1$, then this definition coincides with Lurie's definition of the $\infty$-operad $\mathbb E_1^\otimes$ from \cite[5.1.0.2]{HA}.
\end{example}

\subsection{Algebras and symmetric monoidal functors}
Next we discuss the notion of an algebra and of a symmetric monoidal functor. These are both special cases of a map of $\infty$-multicategories, for which we need the next notion.

\begin{definition}
    A morphism $\alpha\colon I_+ \to J_+$ is \emph{inert} if for each $j \in J$, $\alpha^{-1}(j)$ has exactly one point. 
    Intuitively, these are the permutations and the maps that collapse some points to the basepoint. 
    A morphism is \emph{active} if $\alpha^{-1}(*) = *$. 
    Intuitively, these are the maps that permute, include additional points, and/or fold non-basepoints together. 
\end{definition}

\begin{definition}[{\cite[2.1.2.7]{HA}}]
   Let $\mathbfcal{O}^{\otimes}$ and $\mathbfcal{O}'^{\otimes}$ be $\infty$-multicategories. 
   A \emph{morphism of $\infty$-multicategories} is a map of simplicial sets $\mathbfcal{O}^\otimes \to \mathbfcal{O}'^{\otimes}$ that commutes with the projections to $N(\Fin_*)$ and preserves the cocartesian morphisms over inert morphisms $\alpha\colon I_+ \to J_+$ in $\Fin_*$.
    
    When the map is of the form $\mathbfcal{O}^{\otimes} \to \mathbfcal{C}^{\otimes}$ for a symmetric monoidal $\infty$-category $\mathbfcal{C}^{\otimes}$, we call this an \textit{$\mathbfcal O$-algebra in $\mathbfcal C$}. 
    The set of all $\mathbfcal O$-algebras in $\mathbfcal C$ forms an $\infty$-category $\Alg_{\mathbfcal O}(\mathbfcal C)$. 
    An $n$-simplex in $\Alg_{\mathbfcal O}(\mathbfcal C)$ is a map of simplicial sets $\Delta^n \times \mathbfcal{O}^\otimes \to \mathbfcal{C}^\otimes$ that respects the projection to $N(\Fin_*)$ and at each vertex of $\Delta^n$ is a map of $\infty$-multicategories.

    When the map is of the form $\mathbfcal{C}^{\otimes} \to \mathbfcal{D}^{\otimes}$ for two symmetric monoidal $\infty$-categories $\mathbfcal{C}^{\otimes}$ and $\mathbfcal{D}^{\otimes}$, we call it a \textit{lax symmetric monoidal functor}. 
    The $\infty$-category of such is denoted $\Fun^{\lax}(\mathbfcal{C}^{\otimes},\mathbfcal{D}^{\otimes})$. 
    If the map preserves all of the cocartesian morphisms, not just the inert ones, then it is called a \textit{strong symmetric monoidal functor}, and the $\infty$-category of such is denoted 
    \(
    	\Fun^\otimes(\mathbfcal{C}^{\otimes},\mathbfcal{D}^{\otimes}).
	\)
\end{definition}

\begin{example}
    Let $\Comm^\otimes$ be the $\infty$-category $N(\Fin_*)$ with the identity projection to $N(\Fin_*)$. 
    This arises from \cref{ex:infty_op_from_simplicial_op} applied to the commutative operad: the operad in sets that has a single point at each level, seen as a constant simplicial operad. 
    The $\infty$-category $\Comm^\otimes$ is an $\infty$-multicategory; in fact, it is a symmetric monoidal $\infty$-category.

    A \textit{commutative algebra} in $\mathbfcal{C}$ is {a} $\Comm$-algebra in $\mathbfcal{C}$: a lax symmetric monoidal functor $N(\Fin_*) \rightarrow \mathbfcal{C}^{\otimes}$. 
    On the set $[1]_+$ this functor selects an object $X \in \mathbfcal{C}$, on the fold map $[2]_+ \to [1]_+$ it selects a multiplication $X \otimes X \to X$, and on the map $[0]_+ \to [1]_+$ it selects a unit $\mathbbm{1} \to X$, all of which are well-defined in the homotopy category $h\mathbfcal{C}$. 
    The rest of the functor gives associativity, commutativity, and unitality isomorphisms, and all coherences between these. 
    The $\infty$-category of commutative algebras $\Alg_{\Comm}(\mathbfcal{C})$ is more commonly denoted $\CAlg(\mathbfcal{C})$.
\end{example}

\begin{example}[{\cite[4.1.1.4]{HA}}]
    Let $\Assoc^\otimes$ be the $\infty$-category associated to the associative operad in sets. 
    This is the nerve of the 1-category which has one object for each finite based set $I_+$ and one morphism $I_+ \to J_+$ for each basepoint-preserving function and each linear ordering of each preimage.
	The $\infty$-category $\Assoc^\otimes$ is an $\infty$-multicategory. 
	(In fact, $\Assoc^\otimes$ is a monoidal $\infty$-category \cite[Deﬁnition 4.1.1.10]{HA}, which we do not define here.)
    
    An \textit{(associative) algebra} in $\mathbfcal{C}$ is an $\Assoc$-algebra, in other words a map of $\infty$-multicategories $\Assoc^\otimes \rightarrow \mathbfcal{C}^{\otimes}$.  
    Again, this is a choice of object $X \in \mathbfcal{C}$, multiplication $X \otimes X \to X$ and unit $\mathbbm{1} \to X$ that are well-defined in the homotopy category $h\mathbfcal{C}$, and associativity and unitality isomorphisms and coherences between them. 
    The $\infty$-category of algebras $\Alg_{\Assoc}(\mathbfcal{C})$ is more commonly denoted $\Alg(\mathbfcal{C})$.
\end{example}

By \cite[5.1.0.7]{HA}, there is an equivalence of $\infty$-multicategories $\mathbb E_1^\otimes \cong \Assoc^\otimes$ which hence yields:

\begin{lemma}\label{e1_to_assoc}
    {Let $\mathbfcal{C}$ be a symmetric monoidal $\infty$-category.} There is an equivalence between the $\infty$-categories of algebras
    \[ 
    	\Alg_{\mathbb E_1}(\mathbfcal{C}) \simeq \Alg(\mathbfcal{C}). 
	\]
\end{lemma}

Therefore, algebras in a symmetric monoidal $\infty$-category really are the homotopically meaningful notion of associative algebras.

Let us quickly observe why lax symmetric monoidal functors induce functors on the {$\infty$-}categories of $\mathbfcal{O}$-algebras:

\begin{construction}\label{induced_functor_of_o_algebras}
	For any lax symmetric monoidal functor $F\colon \mathbfcal{C}^{\otimes} \to \mathbfcal{D}^{\otimes}$ and any $\infty$-multicategory $\mathbfcal{O}^\otimes$, the induced map of $\mathbfcal{O}$-algebras is
	\[
		\begin{tikzcd}[column sep = 4em]
			\Fun^{\lax}(\mathbfcal{O}^\otimes,\mathbfcal{C}^{\otimes}) 
				\rar{F \circ -} 
				& 
			\Fun^{\lax}(\mathbfcal{O}^\otimes,\mathbfcal{D}^{\otimes}). 
		\end{tikzcd}
	\]
	This is, by definition, a functor 
	\(
		\Alg_{\mathbfcal{O}}(\mathbfcal{C}) \to \Alg_{\mathbfcal{O}}(\mathbfcal{D})
	\) 
	that we denote $F_{\Alg_{\mathbfcal{O}}}$. 
\end{construction}

We recall from \cite[3.2.4.1 and 3.2.4.4]{HA} a symmetric monoidal structure on  $\CAlg(\mathbfcal{C})$ given by the ``pointwise tensor product.''

\begin{definition}\label{pointwise_tensor}
	We define $\CAlg(\mathbfcal{C})^\otimes$ to be the symmetric monoidal $\infty$-category with the following universal property. 
	Giving a map of simplicial sets $K \to \CAlg(\mathbfcal{C})^\otimes$ over $N(\Fin_*)$ is equivalent to giving a commuting diagram
	\begin{equation}
        \label{universal property of commutative algebras}
		\begin{tikzcd}
		    K \times N(\Fin_*) 
		    	\dar 
				\rar 
				& 
			\mathbfcal{C}^\otimes 
				\dar 
				\\
		    N(\Fin_*) \times N(\Fin_*) 
		    	\rar{\wedge} 
				& 
			N(\Fin_*)
		\end{tikzcd} 
	\end{equation}
	that for each vertex of $K$ over $[m]_+$ sends each cocartesian arrow over an inert morphism $I_+ \to J_+$ to a cocartesian arrow over the resulting inert morphism $([m] \times I)_+ \to ([m] \times J)_+$.
\end{definition}

{To show that this universal property is nonempty, we construct $\CAlg(\mathbfcal{C})^\otimes$ to be the simplicial set whose $n$-simplices are diagrams of the form \cref{universal property of commutative algebras} with $K = \Delta^n$. One can check that this is indeed a symmetric monoidal $\infty$-category.}

This symmetric monoidal structure on $\CAlg(\mathbfcal{C})$ is the cocartesian one by \cite[3.2.4.7]{HA}. An immediate consequence of this is:

\begin{lemma}\label{sm_coproduct_preserving}
    For any {strong} symmetric monoidal functor $\mathbfcal{C}^\otimes \to \mathbfcal{D}^\otimes$, the induced functor on the $\infty$-category of commutative algebras (\cref{induced_functor_of_o_algebras})
    \[ 
    	F_{\CAlg}\colon \CAlg(\mathbfcal{C}) \to \CAlg(\mathbfcal{D}) 
	\]
    lifts to a {strong} symmetric monoidal functor 
    \(
    	F_{\CAlg}\colon \CAlg(\mathbfcal{C})^{\otimes} \to \CAlg(\mathbfcal{D})^{\otimes}.
	\)
\end{lemma}

\section{Transferring algebras from simplicial categories to \texorpdfstring{$\infty$}{∞}-categories}

Let $C_\bullet$ {be} any symmetric monoidal simplicial category. 
In this section we define canonical interchange maps of simplicial sets
\[ 
	N^s(\Alg_{\mathcal{O}_\bullet}(C_\bullet)) \to \Alg_{N^s(\mathcal O_\bullet)}(N^s(C_\bullet)) 
\]
for any simplicial operad $\mathcal{O}_\bullet$. 
In other words, we prove the variant of \cref{thm-alg-interchange-intro} where instead of a model category we have any simplicial category with a symmetric monoidal structure. 
If $C_\bullet$ is fibrant, or if we compose this interchange with the fibrant replacement $C_\bullet \to \Ex^\infty C_\bullet$ from \cref{monoidal_fibrant_replacement}, then this tells us that every algebra in $C_\bullet$ passes to an algebra in the associated $\infty$-category, and moreover the space of maps of algebras goes to the space of maps between them in the $\infty$-category.

These interchange maps are well-known in the case of a symmetric monoidal 1-category \cite[2.1.3.3]{HA}, but we have not been able to find them stated explicitly in the literature for a symmetric monoidal simplicial category. 

\subsection{Spaces of algebra morphisms in a simplicial category}

First recall that when $C_\bullet$ is a symmetric monoidal simplicial $1$-category, the category of algebras $\Alg_{\mathcal{O}_\bullet}(C_\bullet)$ is simplicially enriched in a canonical way. 
This is fairly standard, see e.g. \cite[VII.2.9]{ekmm} or \cite[4.1.8.5]{HA}, but we will spell out the enrichment in detail.

\begin{proposition}\label{enrich_algebras}
    Let $C_\bullet$ be a symmetric monoidal simplicial category. 
    For any simplicial operad $\mathcal{O}_\bullet$, the $1$-category $\Alg_{\mathcal{O}_\bullet}(C_\bullet)$ is simplicially enriched and symmetric monoidal under the tensor product of algebras. 
    In the commutative case, this tensor product is also the coproduct in $\CAlg(C_\bullet)$.
\end{proposition}

{
\begin{remark}
\label{too many bullets}
Although the most honest notation for the simplicial enrichment of the $1$-category $\Alg_{\mathcal{O}_\bullet}(C_\bullet)$ would be $\Alg_{\mathcal{O}_\bullet}(C_\bullet)_\bullet$, we will avoid the proliferation of placeholder indices and simply write $\Alg_{\mathcal{O}_\bullet}(C_\bullet)$ for the symmetric monoidal simplicial category as well as the $1$-category. Note that this is an abuse of notation that departs from the conventions we set to distinguish $1$-categories and simplicial categories.
\end{remark}
}

To be more explicit, for algebras $R$ and $S$ over the operad $\mathcal O_\bullet$, we define the algebra structure on $R \otimes S$ by
\[ 
	\begin{tikzcd}[row sep = 1em]
    	\mathcal{O}_\bullet(n) 
			\rar{(\id,\id)} 
			& 
		\mathcal{O}_\bullet(n) \times \mathcal{O}_\bullet(n) 
			\rar{\mu_R \times \mu_S} 
			\ar[d, phantom, ""{coordinate, name=A}]
			& 
		C_\bullet(R^{\otimes n},R) \times C_\bullet(S^{\otimes n},S)
			\ar[dl, rounded corners, 
				to path={ -- ([xshift=12mm]\tikztostart.east)
					  |- (A)
					  -| ([xshift=-5mm]\tikztotarget.west)
					  -- (\tikztotarget)}]
			\\
			& 
		C_\bullet((R^{\otimes n}) \otimes (S^{\otimes n}),R \otimes S) 
			\rar{\cong} 
			& 
		C_\bullet((R\otimes S)^{\otimes n},R \otimes S).
\end{tikzcd}
\]
In the special case of an associative or commutative algebra, this means the multiplication and unit maps are given by
\[ 
	\begin{tikzcd}[column sep = 4em]
    	(R \otimes S) \otimes (R \otimes S) \cong (R \otimes R) \otimes (S \otimes S) 
			\rar{\mu_R \otimes \mu_S} 
			& 
		R \otimes S,
	\end{tikzcd}
\]
\[
	\begin{tikzcd}[column sep = 4em]
    	\mathbbm{1} \cong \mathbbm{1} \otimes \mathbbm{1} 
	    	\rar{\eta_R \otimes \eta_S} 
			& 
		R \otimes S.
	\end{tikzcd} 
\]
For the enrichment, $\Alg_{\mathcal{O}_\bullet}(C_\bullet)(R,S)$ is taken to be a simplicial subset of $C_\bullet(R,S)$. 
At each simplicial level $k$, it is the set of $f \in C_k(R,S)$ such that both branches of the following diagram agree on $(\rho,f)$, for every $n$ and every $k$-simplex of $n$-fold multiplications $\rho \in \mathcal O_k(n)$:
\[
	\begin{tikzcd}
	    \mathcal O_k(n) \times C_k(R,S) 
	    	\rar{\mu_R \times \id} 
			\dar{\Delta} 
			& 
		C_k(R^{\otimes n},R) \times C_k(R,S) 
			\rar{\circ} 
			& 
		C_k(R^{\otimes n},S) 
			\arrow[d, equals] 
			\\
	    \mathcal O_k(n) \times C_k(R,S)^n 
	    	\rar{\mu_S \times \otimes} 
			& 
		C_k(S^{\otimes n},S) \times C_k(R^{\otimes n},S^{\otimes n}) 
			\rar{\circ} 
			& 
		C_k(R^{\otimes n},S).
	\end{tikzcd}
\]
In the special case of an associative or commutative algebra, this means that the following two diagrams commute:
\[
	\begin{tikzcd}
    	C_k(R,S) 
			\rar{- \circ \mu_R} 
			\dar{\Delta} 
			& 
		C_k(R \otimes R,S) 
			& 
		C_k(R,S) 
			\dar 
			\rar{\circ \eta_R} 
			& 
		C_k(\mathbbm{1},S)
			\arrow[d, equals] 
			\\
	    C_k(R,S) \times C_k(R,S) 
	    	\rar{\otimes} 
			& 
		C_k(R \otimes R,S \otimes S) 
			\uar[swap]{\mu_S \circ -} 
			& 
		{*} 
			\rar{\eta_S} 
			& 
		C_k(\mathbbm{1},S).
	\end{tikzcd}
\]
One might get some intuition for this definition by thinking about it one simplicial level at a time. 
For each value of $k$, $C_k$ is just a symmetric monoidal $1$-category, and $\mathcal O_k(-)$ is an operad in sets. 
As $k$ varies, the {1-}categories $C_k$ have the same objects. 
An $\mathcal O_\bullet$-algebra is an object in this common set of objects, that has the structure of an $\mathcal O_k(-)$-algebra in $C_k$ for all $k$, respecting the face and degeneracy maps. 
At each level $k$ we therefore get a notion of map of algebras in $C_k$, and as $k$ varies these form the simplicial set $\Alg_{\mathcal{O}_\bullet}(C_\bullet)(R,S)$.

\begin{remark}\label{top_case_1}
    If $C_\bullet$ arises by taking singular simplices of a  symmetric monoidal topological $1$-category, then this definition of the mapping space becomes much simpler: it is just the (singular simplices of the) subspace of the topological space $C(R,S)$, consisting of all points $f \in C(R,S)$ that are maps of $|\mathcal O_\bullet|$-algebras, in the sense that both branches of the following diagram agree on $(\rho,f)$ for every $\rho \in |\mathcal O_\bullet(n)|$:
	\[
	    \begin{tikzcd}
    	    {|\mathcal O_\bullet(n)| \times C(R,S)} 
	        	\rar{\mu_R \times \id} 
				\dar{\Delta} 
				& 
			C(R^{\otimes n},R) \times C(R,S) \rar{\circ} 
				& 
			C(R^{\otimes n},S) 
				\arrow[d, equals] 
				\\
	        {|\mathcal O_\bullet(n)| \times C(R,S)^n} 
    	    	\rar{\mu_S \times \otimes} 
				& 
			C(S^{\otimes n},S) \times C(R^{\otimes n},S^{\otimes n}) 
				\rar{\circ} 
				& 
			C(R^{\otimes n},S).
	    \end{tikzcd}
	\]
    In this case, the simplicial mapping sets of $\Alg_{\mathcal{O}_\bullet}(C_\bullet)$ are all Kan complexes,\ ie, $\Alg_{\mathcal{O}_\bullet}(C_\bullet)$ is fibrant.
\end{remark}

We spelled out this structure explicitly so that we could check the following lemma.
\begin{lemma}
    If $R$ and $S$ are two $\mathcal O_\bullet$-algebras in $C_\bullet$, then $\Alg_{\mathcal{O}_\bullet}(C_\bullet)(R,S)$ is isomorphic to the simplicial set in which a $k$-simplex is a functor of 1-categories over $\Fin_*$
    \[ 
    	\{ 0 \to 1 \} \times \mathcal O_k^\otimes \to C_k^\otimes, 
	\]
    whose restrictions to $\{0\} \times \mathcal O_k^\otimes$ and $\{1\} \times \mathcal O_k^\otimes$ are the functors $R^\otimes_k$ and $S^\otimes_k$ that describe the algebra structures on $R$ and $S$ in $C_k$.
\end{lemma}

\begin{proof}
    This is a 1-categorical statement about algebras over an operad valued in sets; the interested reader can check this directly.
\end{proof}

\subsection{The interchange map}

Before we proceed to define the interchange map, we elaborate a bit more on the simplicial nerve $N^s$.  
For each natural number $n$, let $\mathfrak C[\Delta^n]$ be the simplicial category from \cite[1.1.5.1]{HTT}. 
Its objects are $0, \dotsc, n$, and for $0 \leq i,j \leq n$ the mapping space is 
\[
	\textup{Map}_{\mathfrak{C}[\Delta^n]}(i,j) = N(P_{i,j})
\] 
if $i\leq j$, where $P_{i,j}$ is the poset of subsets of $\{i, \dots, j \}$ containing $i$ and $j$, and $\emptyset$ otherwise.

For $C_\bullet$ a simplicial category, $N^s(C_{\bullet})$ is defined to be the simplicial set whose $n$-simplices are the functors of simplicial categories $\mathfrak C[\Delta^n] \to C_\bullet$. 

We can view $\mathfrak{C}[\Delta^{(-)}]$ as a functor  
\[ 
	\mathfrak{C} \colon \Delta \to \Fun(\Delta^{\op},\Cat) = \sCat
\] 
which induces a functor
\[
	\mathfrak{C}[-] \colon \sSet \to \sCat
\]
on the associated category of presheaves. 
The latter functor is left adjoint to the simplicial nerve $N^s$, and this is an example of a nerve-realization adjunction \cite[Section 3]{Kan1958}.

\begin{corollary}
    The simplicial nerve $N^s(\Alg_{\mathcal{O}_\bullet}(C_\bullet))$ has an $n$-simplex for each functor of simplicial categories over $\Fin_*$
    \[ 
    	\mathfrak C[\Delta^n] \times \mathcal O_\bullet^\otimes \to C_\bullet^\otimes, 
	\]
    that on each object of $\mathfrak C[\Delta^n]$ gives a map of multicategories $\mathcal O_\bullet^\otimes \to C_\bullet^\otimes$.
\end{corollary}

\begin{proof}
    This follows by applying the previous lemma to each $k$-simplex in each mapping space in the simplicial category $\mathfrak C[\Delta^n]$.
\end{proof}

\begin{theorem}\label{alg_interchange}
    For each fibrant symmetric monoidal simplicial category $C_\bullet$ and each fibrant simplicial operad $\mathcal O_\bullet$ there is a canonical interchange map
    \[ 
    	\mathcal N\colon N^s(\Alg_{\mathcal{O}_\bullet}(C_\bullet)) 
			\to 
		\Alg_{N^s(\mathcal{O}_\bullet)}(N^s(C_\bullet)) 
	\] 
    that for each functor $\mathfrak C[\Delta^n] \times \mathcal O_\bullet^{{\otimes}} \to C_\bullet^\otimes$ takes the simplicial nerve, and then composes with the unit of the adjunction $\mathfrak{C}[-] \dashv N^s$:
    \[ 
    	\Delta^n \times N^s(\mathcal O_\bullet^\otimes) 
			\to 
		N^s(\mathfrak C[\Delta^n]) \times N^s(\mathcal O_\bullet^\otimes)
    		\cong 
		N^s(\mathfrak C[\Delta^n] \times \mathcal O_\bullet^\otimes) \to N^s(C_\bullet^\otimes). 
	\]
\end{theorem}

This follows directly from the previous results. 
Note that on objects, this map simply represents each algebra by a functor $\mathcal O_\bullet^\otimes \to C_\bullet^\otimes$ and then takes the simplicial nerve. 
On morphisms, the map represents each algebra map as a functor $\mathfrak C[\Delta^1] \times \mathcal O_\bullet^\otimes \to C_\bullet^\otimes$ and then takes the simplicial nerve, giving
\[ 
	\Delta^1 \times N^s(\mathcal O_\bullet^\otimes) \to N^s(C_\bullet^\otimes) 
\]
by the above fact that $N^s$ is a right adjoint and thus preserves products. 
The definition for the simplices of dimension 2 and higher is more complex.

Let $C_\bullet$ be a fibrant symmetric monoidal simplicial category. As a special case of the above theorem, we get a canonical map for algebras over the little intervals operad,
    \[ 
    	\mathcal N\colon N^s(\Alg_{D_1}(C_\bullet)) \to \Alg_{\mathbb E_1}(N^s(C_\bullet)). 
	\] 
This last term is equivalent to $\Alg(N^s(C_\bullet))$ {because $C_\bullet$ is fibrant}.

\begin{corollary}\label{underlying_multiplication}
    Along the canonical interchange map of \cref{alg_interchange}, each $\mathcal O_\bullet$-algebra $X$ in $C_\bullet$ is sent to an $N^s(\mathcal O_\bullet)$-algebra in $N^s(C_\bullet)$ that has the same multiplication map in the homotopy category $hC_\bullet \cong hN^s(C_\bullet)$.
\end{corollary}

\begin{proof}
    Note that the homotopy categories of $C_\bullet$ and $N^s(C_\bullet)$ are canonically identified, since they have the same objects, and the morphisms are equivalence classes of morphisms from $C_0$ by the same relation. 
    The proof of this corollary amounts to writing down definitions: for an algebra $X$ in $C_\bullet$, its multiplication map $\mu\colon X \otimes X \to X$ is a morphism $\mu \in C_0$ that passes to a morphism in the homotopy category. 
    In $N^s(C_\bullet^{\otimes})$, this multiplication becomes a morphism $\{X,X\}_{[2]_+} \to \{X\}_{[1]_+}$ over the basepoint-preserving function $[2]_+ \to [1]_+$ that folds the two points into one. 
    We then factor this into the cocartesian arrow $\{X,X\}_{[2]_+} \to \{X \otimes X\}_{[1]_+}$ and the multiplication map $\mu\colon X \otimes X \to X$. 
    By definition the multiplication in the homotopy category is the second map in this factorization, and is therefore also given by $\mu$.
\end{proof}

We will need to know one more fact about this interchange:

\begin{proposition}\label{interchange_coproucts}
    If the mapping spaces of $C_\bullet$ and $\CAlg(C_\bullet)$ are Kan complexes, then the interchange for commutative algebras
    \[ 
    	\mathcal N\colon N^s(\CAlg(C_\bullet)) \to \CAlg(N^s(C_\bullet)) 
	\]
    preserves coproducts.
\end{proposition}

The assumption about Kan complexes is necessary to ensure that this is a map between $\infty$-categories and not just between simplicial sets. 
Without this, it doesn't even make sense to talk about coproducts!

{Before we prove \cref{interchange_coproucts}, we first recall the 2-functoriality of the simplicial nerve. Let $F, G \colon C_\bullet \to D_\bullet$ be functors between simplicial categories. Recall that a natural transformation between $F$ and $G$ is equivalent to a functor $C_\bullet \times \mathbf{1} \to D_\bullet$ that restricts to $F$ on $C_\bullet \times \{0\}$ and $G$ on $C_\bullet \times \{1\}$, where $\mathbf{1}$ is the one-arrow $1$-category. Since $N^s(\mathbf{1}) \cong N(\mathbf{1}) \cong \Delta^1$ and $N^s$ preserves products, when we apply $N^s$ this gives a morphism of simplicial sets 
\[
	N^s(C_\bullet) \times \Delta^1 \to N^s(D_\bullet) 
\]
that restricts to $N^s(F)$ and $N^s(G)$ on the endpoints of the interval $\Delta^1$. When $C_\bullet$ and $D_\bullet$ are fibrant, then this is exactly a natural transformation between the $\infty$-categories $N^s(C_\bullet)$ and $N^s(D_\bullet)$. Hence $N^s$ sends natural transformations between fibrant simplicial categories to natural transformations of $\infty$-categories.}

\begin{proof}[{Proof of \cref{interchange_coproucts}}]
	Let $R$ and $S$ be commutative algebras in $C_\bullet$, represented by functors
	\(
		R^\otimes\colon \Fin_* \to C_\bullet^\otimes
	\)
	and
	\(
		S^\otimes\colon \Fin_* \to C_\bullet^\otimes. 
	\)
    The coproduct is $R \otimes S$, with multiplication given in the usual way by shuffling. This defines a functor
    \[ 
    	(R \otimes S)^\otimes\colon \Fin_* \to C_\bullet^\otimes. 
	\]
    Applying the simplicial nerve, we get a functor
    \[ 
    	N^s((R \otimes S)^\otimes)\colon N(\Fin_*) \to N^s(C_\bullet^\otimes) 
	\]
    that receives natural transformations over $N(\Fin_*)$ from the two functors
    \[ 
    	N^s(R^\otimes)\colon N(\Fin_*) \to N^s(C_\bullet^\otimes),
	    	\qquad 
	    N^s(S^\otimes)\colon N(\Fin_*) \to N^s(C_\bullet^\otimes). 
	\]
    By \cite[3.2.4.7]{HA}, the symmetric monoidal structure on $\CAlg(N^s(C_\bullet))$ is the cocartesian one, so it suffices to show $N^s(R^\otimes) \otimes N^s(S^\otimes) \simeq N^s((R \otimes S)^\otimes)$. 
    The left hand side is constructed as follows. 
    We pull back the functors $N^s(R^\otimes)$ and $N^s(S^\otimes)$ along the two inert morphisms $[2]_+ \to [1]_+$ to get a single functor
    \[ 
    	(N^s(R^\otimes), N^s(S^\otimes))\colon [2]_+ \times N(\Fin_*) \to N^s(C_\bullet^\otimes) 
	\]
    lying over the map $[2]_+ \times N(\Fin_*) \to N(\Fin_*)$ sending $([2]_+,[m]_+) \mapsto ([2m]_+)$.
    This is the formal way of pairing our two functors. 
    Then, as also discussed below \cref{df:sm_infinity_cat}, we push forward along the active morphism $[2]_+ \to [1]_+$ to get a single functor
    \[ 
    	N^s(R^\otimes) \otimes N^s(S^\otimes)\colon [1]_+ \times N(\Fin_*) \to N^s(C_\bullet^\otimes) 
	\]
    representing the symmetric monoidal product, and hence in this case the coproduct (on the level of the fiber over $[1]_+$).

    However, we can do each of these steps on the underlying simplicial category and then take the simplicial nerve. 
    Consider the functor
    \[ 
    	(R,S)^\otimes\colon [2]_+ \times \Fin_* \to C_\bullet^\otimes 
	\]
    that takes every finite set $[m]_+$ to a list of $2m$ objects, consisting of $m$ copies of $R$ and $m$ copies of $S$. 
    Each map $[m]_+ \to [n]_+$ applies the corresponding multiplications to these copies of $R$ and $S$ separately. 
    This is indeed a lift of the two functors $R^\otimes$ and $S^\otimes$ along the two inert morphisms $[2]_+ \to [1]_+$ at the level of simplicial categories, so it is also true after we take the simplicial nerve.
    
    Similarly, the functor
    \[ 
    	(R \otimes S)^\otimes\colon [1]_+ \times \Fin_* \to C_\bullet^\otimes
	\]
    is a pushforward of $(R,S)^\otimes$ along the active map $[2]_+ \to [1]_+$ at the level of simplicial categories, so it is also true after we take the simplicial nerve. 
    So we do indeed get the same result.

    This entire comparison is natural in $R$ and $S$, so if we apply the unit map $\mathbbm{1} \to R$ in the place of $R$, this shows that the inclusion $S^\otimes \to (R \otimes S)^\otimes$ goes to the inclusion $N^s(S^\otimes) \to N^s((R \otimes S)^\otimes)$, and similarly for $R$. 
    Therefore everything is compatible, so the interchange map preserves coproducts.
\end{proof}

\section{Transferring algebras from model categories to \texorpdfstring{$\infty$}{∞}-categories}

To finish the proof of \cref{thm-alg-interchange-intro}, it remains to explain how the underlying symmetric monoidal $\infty$-category of a symmetric monoidal simplicial model category $M_\bullet$ is related to the simplicial nerve $N^s(M_\bullet)$. 
Recall from \cref{underlying_sm_inf_cat_1} that the underlying symmetric monoidal $\infty$-category is usually defined as {the} simplicial nerve of the simplicial category of bifibrant objects $N^s((M^{\mathrm{cf}}_\bullet)^\otimes)$. 
However, in this section we will give several more equivalent definitions of the underlying symmetric monoidal $\infty$-category, including one that does not even mention the cofibrant and fibrant objects. The key to building these is to use Dwyer--Kan localization.

\subsection{Dwyer--Kan Localization}
Recall that a \emph{Dwyer--Kan localization}, or quasi-localization, is a functor $\mathbfcal C \to \mathbfcal C[W^{-1}]$ formed by selecting an $\infty$-category $\mathbfcal C$ and a class of morphisms $W \subseteq \mathbfcal C_1$ to invert, and then constructing the universal $\infty$-category receiving a map from $\mathbfcal C$ in which each morphism in $W$ is sent to an isomorphism. 
We use the same notation $C_\bullet \to C_\bullet[W^{-1}]$ for the Dwyer--Kan localization of a simplicial category, which can be constructed directly via the hammock localization \cite{dk3}. 
On simplicial nerves, the Dwyer--Kan localization of simplicial categories satisfies the same universal property as the localization in $\infty$-categories:
\[ 
	\begin{tikzcd}[row sep = 0em]
    		& 
		N^s(C_\bullet)[W^{-1}] 
			\ar{dd}{\sim} 
			\\
	    N^s(C_\bullet) 
	    	\ar[in=180]{ru} \ar[out=-45,in=180]{rd} 
			& 
			\\
	    	& 
	    N^s(C_\bullet[W^{-1}])
	\end{tikzcd} 
\]
Note that Dwyer--Kan localization in this sense is more general and less structured than Bousfield localization, also called reflective localization, and which in \cite{HTT,HA} is simply called ``localization.'' 
We will not need to use Bousfield localization in this paper; when we write ``localization,'' we are always referring to the Dwyer--Kan localization.

In general, Dwyer--Kan localization does not preserve symmetric monoidal structure. 
However, we have the following consequence of \cite[Thm 3.2.2]{hinich} as stated in \cite[Prop A.5]{nikolaus_scholze}:
\begin{proposition}\label{sm_localization}
    Suppose $\mathbfcal{C}^\otimes$ is {a} symmetric monoidal $\infty$-category, and $W$ is a class of morphisms in the underlying $\infty$-category $\mathbfcal{C}$, with the property that tuples of morphisms in $W$ are preserved by the pushforward functors $\alpha_!$ for every active morphism $\alpha\colon I_+ \to J_+$. 
    Then there is a canonical ``strict symmetric monoidal localization'' $\mathbfcal{C}^\otimes \to \mathbfcal{C}^\otimes[(W^\otimes)^{-1}]$, that is the universal {strong} symmetric monoidal functor from $\mathbfcal{C}^\otimes$ that sends each map in $W$ to an isomorphism. 
    Furthermore, on the underlying $\infty$-category, it has the universal property of the Dwyer--Kan localization $\mathbfcal C \to \mathbfcal C[W^{-1}]$.
\end{proposition}

To be more specific, $W^\otimes$ is a collection of maps in $\mathbfcal{C}^\otimes$ that all lie over identity morphisms in $\Fin_*$. 
In other words, these maps all live in the various fiber categories of $\mathbfcal{C}^\otimes$, and in each fiber category $\mathbfcal{C}^\otimes_{{[}n{]}_+} \simeq \mathbfcal{C}^{\times n}$ it is just the collection $W^{\times n}$ of $n$-tuples of maps in $W$. 
So the symmetric monoidal localization is, in each fiber category, equivalent to the $n$-fold product of the Dwyer--Kan localizations $\mathbfcal{C} \to \mathbfcal{C}[W^{-1}]$. 
This motivates the following more general construction.

\begin{construction}\label{sm_localization_attempt}
    For any symmetric monoidal $\infty$-category $\mathbfcal{C}^\otimes$ with a class of weak equivalences $W \subseteq \mathbfcal C_1$, we define $\mathbfcal{C}^\otimes[(W^\otimes)^{-1}]$ to be any Dwyer--Kan localization that universally inverts the arrows in $W^\otimes$, along with a map to $N(\Fin_*)$ such that the following diagram commutes.
    \[
    	\begin{tikzcd}
        	\mathbfcal{C}^\otimes 
				\dar 
				\rar 
				& 
			\mathbfcal{C}^\otimes[(W^\otimes)^{-1}] 
				\dar 
				\\
	        N(\Fin_*) 
	        	\arrow[r,equals] 
				& 
			N(\Fin_*)
    	\end{tikzcd} 
    \]
    This always exists by the universal property of Dwyer--Kan localization, because every arrow in $W^\otimes$ goes to an isomorphism in $N(\Fin_*)$. 
    Note that by \cite[2.3.1.5]{HTT}, the functor from $\mathbfcal{C}^\otimes[(W^\otimes)^{-1}]$ to $N(\Fin_*)$ is automatically an inner fibration.
\end{construction}

One might na\" ively think that the above construction allows us to take any symmetric monoidal $\infty$-category, invert any class of equivalences, and get another symmetric monoidal $\infty$-category. 
It does not. 
The issue, and the reason why \cref{sm_localization} is a nontrivial result, is that when we invert the maps in $W^\otimes$, in general we invert \emph{more} in each fiber than just the maps in $W^{\times n}$. 
This means that the resulting map $\mathbfcal{C}^\otimes[(W^\otimes)^{-1}] \to N(\Fin_*)$ no longer has the correct fibers, preventing it from being a cocartesian fibration with the appropriate Segal condition.

However, we will show in \cref{underlying_sm_inf_cat_compare_2} that \cref{sm_localization_attempt} can still produce a symmetric monoidal $\infty$-category in some cases where the assumption of \cref{sm_localization} is not fulfilled. 
The following two lemmas are useful for this.

\begin{lemma}\label{sm_equivalences}
    If $\mathbfcal{C}^\otimes$ is a symmetric monoidal $\infty$-category, $\mathbfcal{D}$ is any $\infty$-category with a functor to $N(\Fin_*)$, and we have an equivalence $\mathbfcal{C}^\otimes \to \mathbfcal{D}$ or $\mathbfcal{D} \to \mathbfcal{C}^\otimes$ over $N(\Fin_*)$, then $\mathbfcal{D}$ is a symmetric monoidal $\infty$-category and the equivalence with $\mathbfcal{C}^\otimes$ is a {strong} symmetric monoidal functor.
\end{lemma}

\begin{proof}
    Since $\mathbfcal{D}$ is an $\infty$-category, the functor $\mathbfcal{D} \to N(\Fin_*)$ is automatically an inner fibration by \cite[2.3.1.5]{HTT}. 
    The criterion for whether an arrow in $\mathbfcal{D}$ is cocartesian is expressed in terms of the homotopy type of its mapping spaces by \cite[2.4.4.3]{HTT}. 
    As these are unchanged under equivalence of $\infty$-categories, $\mathbfcal{D} \to N(\Fin_*)$ is a cocartesian fibration as well, and the map to or from $\mathbfcal{C}^{\otimes}$ preserves these cocartesian arrows, so by definition it is a {strong} symmetric monoidal functor.
\end{proof}

\begin{lemma}\label{equivalence_fiberwise}
    A map of symmetric monoidal $\infty$-categories $\mathbfcal{C}^\otimes \to \mathbfcal{D}^\otimes$ is an equivalence of $\infty$-categories (after forgetting the map to $N(\Fin_*)$) iff it is an equivalence on the underlying $\infty$-categories $\mathbfcal{C} \simeq \mathbfcal{D}$.
\end{lemma}

\begin{proof}
    This follows quickly from the fact that equivalences of cocartesian fibrations are detected fiberwise \cite[\href{https://kerodon.net/tag/023M}{Theorem 023M}]{kerodon}{, and the fact that there is a canonical equivalence $\mathbfcal{C}^\otimes_{[n]_+}\simeq (\mathbfcal{C}^{\otimes}_{[1]_+})^{\times n}$.}
\end{proof}

\subsection{Alternate models for the underlying symmetric monoidal \texorpdfstring{$\infty$}{∞}-category}

For the non-formal input we will use the following result, which has been proven many times but which originally comes from Dwyer and Kan. {Recall from \cite{bergner2007} that a \emph{Dwyer--Kan equivalence} or \emph{DK-equivalence} of simplicial categories is a simplicial functor $F\colon C_\bullet \to D_\bullet$ that induces weak homotopy equivalences of simplicial sets $C_\bullet(x,y) \to D_\bullet(F(x),F(y))$ and whose map of homotopy categories $\pi_0 F\colon \pi_0 C_\bullet \to \pi_0 D_\bullet$ is essentially surjective.}

\begin{proposition}[{\cite[Propositions 4.8 and 5.2]{dk3}}]
\label{dk_underlying_inf_cat}
    For any simplicial model category $M_\bullet$ with underlying model category $M_0$, if $M^\flat_\bullet$ is any full subcategory containing the cofibrant objects $M^c_\bullet$, then the maps of simplicial categories
    \[ 
    	\begin{tikzcd}[cramped]
	    	M^{\mathrm{cf}}_\bullet 
				\ar[r] 
				&
			M^\flat_\bullet 
				\ar[r]
				&
			M_\bullet 
				&
			M_0 
				\ar[l]
		\end{tikzcd}
	\]
    induce Dwyer--Kan equivalences after localizing
    \[ 
    	\begin{tikzcd}[cramped]
			M^{\mathrm{cf}}_\bullet 
				\ar[r, "\sim"]
				&
			M^\flat_\bullet[(W^\flat)^{-1}] 
				\ar[r, "\sim"]
				&
			M_\bullet[W^{-1}] 
				&
			M_0[W^{-1}]. 
				\ar[l, "\sim"']
		\end{tikzcd}
	\]
    Here $W^\flat = W \cap M^\flat_0$ is just the weak equivalences restricted to the subcategory $M^\flat_0$.
\end{proposition}

{Our preference for cofibrant (as opposed to fibrant) objects here will be explained later in \cref{warning cofibrant bias}.}

{In light of \cref{dk_underlying_inf_cat}}, if all of the mapping spaces $M_\bullet(a,b)$ in a simplicial model category $M_\bullet$ are Kan complexes, or we are willing to apply $\Ex^\infty$ to them first, then taking the simplicial nerve gives equivalences of $\infty$-categories
\[ 
	N^s(M^{\mathrm{cf}}_\bullet) 
		\simeq 
	N^s(M^\flat_\bullet)[(W^\flat)^{-1}] 
		\simeq 
	N^s(M_\bullet)[W^{-1}] 
		\simeq 
	N(M_0)[W^{-1}]. 
\]
{(As in \cref{def:infty_category}, this follows from the fact that $N^s$ is a right Quillen functor \cite[Section 1.3]{Joyal2007} for the Bergner model structure on simplicial categories from \cite{bergner2007}, in which the weak equivalences are the Dwyer--Kan equivalences.)}

Now assume that
\begin{itemize}
    \item $M_\bullet$ is a simplicial model category and a symmetric monoidal simplicial category, 
    \item $M^\flat_\bullet$ is a full subcategory containing the cofibrant objects, and
    \item the tensor product preserves objects in $M^\flat_\bullet$ and weak equivalences between them.
\end{itemize}

Note this is weaker than the assumption that $M_\bullet$ is a symmetric monoidal simplicial model category, because we do not assume that the model structure is compatible with the symmetric monoidal structure in the usual way. Instead, we only assume that the tensor preserves the ``weakly cofibrant'' objects $M^\flat_\bullet$ and the weak equivalences between them. 

\begin{proposition}\label{underlying_sm_inf_cat_compare}
    Under these assumptions, the localization in \cref{sm_localization_attempt} successfully creates a symmetric monoidal $\infty$-category out of $N^s((\Ex^\infty M^\flat_\bullet)^{\otimes})$, and we get the following commuting diagram of symmetric monoidal $\infty$-categories and {strong} symmetric monoidal functors.
    \[
	    \begin{tikzcd}
        	N^s((M^{\mathrm{cf}}_\bullet)^{\otimes}) 
				\dar{\sim} 
				\rar 
				&
	        N^s((\Ex^\infty M^\flat_\bullet)^{\otimes}) 
	        	\dar 
				\\
	        N^s((M^{\mathrm{cf}}_\bullet)^{\otimes})[((W^{\mathrm{cf}})^\otimes)^{-1}] 
	        	\rar{\sim} 
				& 
			N^s((\Ex^\infty M^\flat_\bullet)^{\otimes})[((W^\flat)^\otimes)^{-1}]
    	\end{tikzcd} 
	\]
\end{proposition}

\begin{proof}
    The categories in the bottom row are symmetric monoidal $\infty$-categories by \cref{sm_localization}, since on the subcategories $M^{\mathrm{cf}}_0$ and $M^\flat_0$ the tensor product preserves all weak equivalences. 
    The left-hand vertical arrow is an equivalence of {$\infty$-categories. To see this, recall that between bifibrant objects in a simplicial model category, every weak equivalence is a homotopy equivalence and therefore becomes an isomorphism in the $\infty$-category defined by the simplicial nerve. Furthermore, Dwyer--Kan localization along a set of isomorphisms produces an equivalent $\infty$-category. Therefore the map of fibers or underlying $\infty$-categories 
    \[
        N^s(M^{\mathrm{cf}}_\bullet)  
        \to N^s(M^{\mathrm{cf}}_\bullet)[(W^{\mathrm{cf}})^{-1}] 
    \]is an equivalence, so by \cref{equivalence_fiberwise} we get an equivalence on the total $\infty$-categories as well. On the other hand,} the right-hand vertical arrow is not an equivalence.
    
    To show that the map along the bottom is a {strong} symmetric monoidal functor, by \cref{sm_equivalences} it suffices to show it is an equivalence. 
    By \cref{equivalence_fiberwise}, it suffices to check this on the underlying $\infty$-categories, where it follows from \cref{dk_underlying_inf_cat}.
\end{proof}

\begin{proposition}\label{underlying_sm_inf_cat_compare_2}
    Under the same assumptions as \cref{underlying_sm_inf_cat_compare}, the localization defined in \cref{sm_localization_attempt} successfully creates a symmetric monoidal $\infty$-category out of both $N^s(M_\bullet^{\otimes})$ and $N(M_0^{\otimes})$, and we have the following diagram of symmetric monoidal $\infty$-categories and symmetric monoidal functors, in which {the solid arrows are strong symmetric monoidal and} the dashed arrows are only lax symmetric monoidal:
    \[
    \begin{tikzcd}
        N^s((\Ex^\infty M^\flat_\bullet)^{\otimes}) 
        	\dar 
			\rar 
			&
        N^s((\Ex^\infty M_\bullet)^{\otimes}) 
        	\dar[dashed] 
			\rar[<-] 
			&
        N(M_0^{\otimes}) 
        	\dar[dashed] 
			\\
        N^s((\Ex^\infty M^\flat_\bullet)^{\otimes})[((W^\flat)^\otimes)^{-1}] 
        	\rar{\sim} 
			&
        N^s((\Ex^\infty M_\bullet)^{\otimes})[(W^\otimes)^{-1}] 
        	\rar[<-]{\sim} 
			&
        N(M_0^{\otimes})[(W^\otimes)^{-1}].
    \end{tikzcd} \]
\end{proposition}

\begin{proof}
    The argument is the same as before, but in place of \cref{sm_localization} we use \cite[Proposition A.14]{nikolaus_scholze}, which tells us that the vertical arrows are lax symmetric monoidal functors that are Dwyer--Kan localizations in both the absolute and fiberwise sense. 
    This in turn allows us to see that the functors along the bottom are also equivalences by \cref{dk_underlying_inf_cat}. 
    The only work here is to verify that the assumptions of \cite[Proposition A.14]{nikolaus_scholze} apply, in other words that the cocartesian fibrations in the top row are left-derivable in the sense of \cite[Definition A.8]{nikolaus_scholze}. 
    For $N(M_0^{\otimes})$ this is already checked directly in \cite[Examples A.10--A.13]{nikolaus_scholze}. 
    For $N^s((\Ex^\infty M_\bullet)^{\otimes})$, the arguments in \cite[Examples A.10--A.13]{nikolaus_scholze} apply equally well, using the fact that simplicial model categories have simplicially enriched cofibrant replacement functors \cite[13.2.4]{riehl_book2}. For $N^s((\Ex^\infty M^\flat_\bullet)^{\otimes})$, this follows from the previous proposition.
\end{proof}

\begin{remark}\label{qx_qy_explanation}
    The difference between the strong and lax symmetric monoidal functors in \cref{underlying_sm_inf_cat_compare_2} can be explained as follows. The tensor products in the $\infty$-categories along the top row are all given by the tensor product $X \otimes Y$ on the underlying model category, while the tensor products along the bottom row are given by the derived tensor product $QX \otimes QY$, where $Q$ denotes cofibrant replacement. The lax symmetric monoidal structure on the vertical functors corresponds to the canonical map $QX \otimes QY \to X \otimes Y$. The vertical functors are not strong symmetric monoidal because this map is not an isomorphism in the $\infty$-category (ie~it is not a weak equivalence in the model category $M$). Only on the far left of the diagram do we assume that $X$ and $Y$ are in the subcategory $M^\flat$, which makes $QX \otimes QY \to X \otimes Y$ a weak equivalence and therefore makes the functor strong symmetric monoidal.
\end{remark}

As a result, we now have three more valid definitions for the underlying symmetric monoidal $\infty$-category:

\begin{definition}\label{underlying_sm_inf_cat_2}
    For any simplicial category $M_\bullet$ with the above assumptions, we also refer to any of the categories in the bottom row of \cref{underlying_sm_inf_cat_compare_2} as the underlying symmetric monoidal $\infty$-category of $M_\bullet$. 
    When the mapping spaces of $M_\bullet$ are Kan complexes, we also extend this honor to the corresponding categories in which the $\Ex^\infty$ has been removed.
\end{definition}

By \cref{underlying_sm_inf_cat_compare} and \cref{underlying_sm_inf_cat_compare_2}, these are equivalent to Lurie's notion of the underlying symmetric monoidal $\infty$-category 
\(
	N^s((M^{\mathrm{cf}}_\bullet)^{\otimes})
\)
from \cref{underlying_sm_inf_cat_1}. 
In practice they tend to be more useful, because we don't have to take cofibrant and fibrant replacements of our objects.

\begin{warning}
	\label{warning cofibrant bias}
    \cref{underlying_sm_inf_cat_compare} also holds when we select a subcategory $M^\sharp_\bullet$ containing the \emph{fibrant} objects, which is closed under the tensor product, and on which the tensor preserves weak equivalences. 
    However, \cref{underlying_sm_inf_cat_compare_2} does not hold in general in this right-derived case. 
    The problem is that it would make the map $N^s(M_\bullet) \to N^s(M_\bullet)[W^{-1}]$ into a lax symmetric monoidal functor. But in the right-derived case this functor ought to be \emph{oplax}, because the natural map goes from the tensor product to the derived tensor product {($X \otimes Y \to RX \otimes RY$)}, not the other way around! (Indeed, we can show that this functor is oplax by applying \cref{underlying_sm_inf_cat_compare_2} to the opposite model category $M^{\op}$.)
\end{warning}

\subsection{The interchange map for model categories}

We may now prove \cref{thm-alg-interchange-intro}. For convenience, we restate the theorem below. 

\begin{theorem}
\label{body_theorem_A}
 Let $M_\bullet$ be a symmetric monoidal simplicial model category with underlying symmetric monoidal $\infty$-category $\mathbfcal{M}$, and let $\mathcal{O}_\bullet$ be a fibrant simplicial operad. Then there is a canonical ``interchange'' map of simplicial sets
    \[ 
    	N^s(\Alg_{\mathcal{O}_\bullet}(M_\bullet)) \to \Alg_{N^s(\mathcal{O}_\bullet)}(\mathbfcal{M}). 
	\]
\end{theorem}

\begin{proof}
For any symmetric monoidal simplicial model category $M$, we form the map of simplicial sets
    \[ 
    	\begin{array}{ll}
            N^s(\Alg_{\mathcal{O}_\bullet}(M_\bullet)) 
            	\to N^s(\Alg_{\mathcal{O}_\bullet}(\Ex^\infty M_\bullet))
				& 
			\textup{\cref{monoidal_fibrant_replacement}, \cref{induced_functor_of_o_algebras}}
				\\[0.3em]
            \phantom{N^s(\Alg_{\mathcal{O}_\bullet}(M_\bullet))} 
            	\to  \Alg_{N^s(\mathcal{O}_\bullet)}(N^s(\Ex^\infty M_\bullet)) 
				& 
			\textup{\cref{alg_interchange}} 
				\\[0.3em]
            \phantom{N^s(\Alg_{\mathcal{O}_\bullet}(M_\bullet))}
            	\to \Alg_{N^s(\mathcal{O}_\bullet)}(N^s(\Ex^\infty M_\bullet)[W^{-1}]) 
				&  
			\textup{\cref{underlying_sm_inf_cat_compare_2}, \cref{induced_functor_of_o_algebras}.}
	    \end{array} 
	\]
The fact that $N^s((\Ex^\infty M_\bullet)^{\otimes})[(W^\otimes)^{-1}]$ is the underlying symmetric monoidal $\infty$-category of $M_\bullet$ finishes the proof.
\end{proof}

\section{Transferring bialgebras and Hopf algebras}

In this section we extend the above definitions and results to commutative bialgebras and Hopf algebras, and prove \cref{main_comparison-intro}.

\subsection{An alternative definition of Hopf algebras}

We begin in the setting of a symmetric monoidal $1$-category $(C, \otimes, \mathbbm{1})$. 
We give {an} alternative characterization of Hopf algebras that will be easier to recreate in simplicial categories and $\infty$-categories. 
Recall first the standard definition of a Hopf algebra: 

\begin{definition}
    A bialgebra $R$ in a symmetric monoidal 1-category $(C,\otimes,\mathbbm{1})$ is an object $R$ that has both the structure of an algebra and a coalgebra, such that the comultiplication $\delta \colon R \to R \otimes R$ and counit $\varepsilon \colon R \to \mathbbm{1}$ are algebra homomorphisms, or equivalently the multiplication $\mu \colon R \otimes R \to R$ and unit $\eta \colon \mathbbm{1} \to R$ are coalgebra homomorphisms.
\end{definition}

\begin{definition}\label{def-hopf}
	Let $R$ be a bialgebra in $(C, \otimes, \mathbbm{1})$. 
	We say that $R$ is a \emph{Hopf algebra} if it has an \emph{antipode} morphism $\alpha \colon R \to R$ such that the following diagram commutes.
	\begin{equation}\label{eq:antipode_condition}
		\begin{tikzcd}
    			& 
		    R \otimes R 
    			\arrow{rr}{\alpha \otimes \id} 
				&  
				&
			R \otimes R 
				\arrow[dr, "\mu"]
				&
				\\
		    R 
		    	\arrow[rr, "\varepsilon"] 
				\arrow[ru, "\delta"] 
				\arrow[rd, "\delta"'] 
				&  
				& 
			\mathbbm{1} \arrow[rr, "\eta"] 
				& 
				& 
			R 
				\\
		    	& 
			R \otimes R 
				\arrow{rr}{\id \otimes \alpha} 
				&  
				&
			R \otimes R 
				\arrow[ur, "\mu"']
		\end{tikzcd}   
	\end{equation}
\end{definition}

It turns out that being a Hopf algebra is a property of a bialgebra, not structure. 
To see this, consider the \emph{shear morphisms}:
\begin{align}
    \begin{tikzcd}[column sep = 4em, ampersand replacement=\&] 
    	\sh \colon R \otimes R 
			\rar{\delta \otimes \id} 
			\& 
		R \otimes R \otimes R
			\rar{\id \otimes \mu} 
			\& 
		R \otimes R. 
	\end{tikzcd} 
	\label{eq:shear_definition}
		\\
    \begin{tikzcd}[column sep = 4em, ampersand replacement=\&] 
    	\phantom{sh \colon} R \otimes R
			\rar{\id \otimes \delta} 
			\& 
		R \otimes R \otimes R
			\rar{\mu \otimes \id} 
			\& 
		R \otimes R. 
	\end{tikzcd}
	\label{eq:left_shear_definition}
\end{align}
We refer to \cref{eq:shear_definition} as the \emph{right shear map} and to \cref{eq:left_shear_definition} as the \emph{left shear map}. 
Invertibility of these maps characterizes those bialgebras in $C$ that are Hopf algebras. 

\begin{proposition}[{\cite[Lemma 5.1]{BLV2011}}]
\label{shear_equivalences}
    Let $R$ be a bialgebra in a symmetric monoidal 1-category.
    The following are equivalent: 
    \vspace*{-1em}
    \begin{enumerate}[label=(\alph*),font=\upshape]
    	\item $R$ is a Hopf algebra as in \cref{def-hopf}.
        \item The right shear map \cref{eq:shear_definition} is an isomorphism.
        \item The left shear map \cref{eq:left_shear_definition}
        is an isomorphism.
    \end{enumerate}
\end{proposition}

Given this proposition, we are therefore justified in taking ``Hopf algebra'' to mean a bialgebra whose shear map is an isomorphism. 
In the homotopical setting of a simplicial category or an $\infty$-category, we will ask for the shear map to be an equivalence, ie\ an isomorphism in the homotopy category. 

{\begin{remark}
	The shear maps above also appear in the context of group-like objects in a cartesian symmetric monoidal category. In a group written additively, the right shear map may be informally described by $(x,y) \mapsto (x, x + y)$, which also describes a shear transformation of the coordinate plane. This motivates the terminology we have chosen.
	In the Hopf algebra literature, the shear maps defined above are often called \emph{fusion operators} \cite[Section 5.2]{BLV2011}. 
\end{remark}}

\subsection{Bialgebras and Hopf algebras in simplicial categories}

Let $C_\bullet$ be a symmetric monoidal simplicial category. {Recall that a simplicial category is a category enriched in simplicial sets; in this description, the opposite category of $C_\bullet$ is the category $C_\bullet^{\op}$ with $C_\bullet^{\op}(a,b) = C_\bullet^{\op}(b,a)$. In the formulation of a simplicial category as a simplicial object in $1$-categories which is constant on objects (\cref{simplicial categories as sCat}), the opposite category may be understood as an instance of the opposite of a simplicial set as described in \cite[Section 1.2.1]{HTT}. }

A \emph{coalgebra} in $C_\bullet$ is an algebra in the opposite category. 
The simplicial category of coalgebras is defined by taking the opposite category twice:
\[ 
	\CoAlg(C_\bullet) := \Alg(C_\bullet^{\op})^{\op}. 
\]
A commutative bialgebra should be an object of $C_\bullet$ that is a commutative algebra and a coalgebra in a compatible way. 
To define this precisely, recall from \cref{enrich_algebras} that $\CAlg(C_\bullet)$ is a symmetric monoidal simplicial category as well, and that the tensor product is the coproduct. 
We deduce that the opposite category $\CAlg(C_\bullet)^{\op}$ is simplicial and Cartesian monoidal, allowing us to define a simplicial category of algebra objects there:
\[ 
	\CBiAlg(C_\bullet) 
		:= \BiAlg_{(\Comm,\Assoc)}(C_\bullet) 
		= \Alg(\CAlg(C_\bullet)^{\op})^{\op}. 
\]
More generally, for any simplicial operad $\mathcal O_\bullet$, we have the following definition.

\begin{definition}
    In any symmetric monoidal simplicial category $C_\bullet$, the simplicial category of $(\Comm,\mathcal{O}_\bullet)$-bialgebras is defined as
    \[ 
    	\BiAlg_{(\Comm,\mathcal{O}_\bullet)}(C_\bullet) 
			:= \Alg_{\mathcal{O}_\bullet}(\CAlg(C_\bullet)^{\op})^{\op}. 
	\]
\end{definition}

We specialize to the case that $\mathcal{O}_\bullet$ is the singular simplices of the little intervals operad, which we denote $D_1$, for the rest of the section. 
As in \cref{rem-comm-D1}, this restriction is not essential for the proofs that follow.

Note that each $(\Comm,D_1)$-bialgebra $X$ in $C_\bullet$ produces multiplication $\mu$ and comultiplication $\delta$ maps in the homotopy category $hC_\bullet$. 
We can compose these to produce right and left shear maps as in \cref{eq:shear_definition,eq:left_shear_definition}, which are morphisms in $hC_\bullet$. 
We say that a bialgebra is a \emph{Hopf algebra} if its shear map is an isomorphism in the homotopy category. 
This non-standard definition is justified by \cref{shear_equivalences}. 

\begin{definition}
	For any symmetric monoidal simplicial category $C_\bullet$, the \emph{simplicial category of commutative Hopf algebras} is the full subcategory
\[ 
	\CHopf(C_\bullet) 
		\subseteq 	
	\Alg(\CAlg(C_\bullet)^{\op})^{\op} 
\]
	spanned by those bialgebras for which the shear map \cref{eq:shear_definition} is an isomorphism in $hC_\bullet$.
\end{definition}

\subsection{Bialgebras and Hopf algebras in \texorpdfstring{$\infty$}{∞}-categories}

We next present the ``abstract'' definition of a bialgebra and a Hopf algebra, following \cite{HopfInfty}.

Let $\mathbfcal{C}^\otimes$ be any symmetric monoidal $\infty$-category. 
As in the concrete case, the category $\CAlg(\mathbfcal{C})$ has all finite coproducts, as discussed in \cref{pointwise_tensor}. 
Along the forgetful functor to the homotopy category $h\mathbfcal{C}$, these agree with the tensor product of commutative algebras that we defined above. 
Since $\CAlg(\mathbfcal{C})$ has all finite coproducts, the opposite $\infty$-category $\CAlg(\mathbfcal{C})^{\op}$ has all finite products, and therefore extends to a Cartesian symmetric monoidal $\infty$-category $(\CAlg(\mathbfcal{C})^{\op})^\times$ in an essentially unique way \cite[2.4.1.9]{HA}. 
Furthermore, product-preserving functors are identified with {strong} symmetric monoidal functors \cite[2.4.1.8]{HA}. 
Since this tensor product is given on the underlying objects by the pointwise tensor product, we will write it as $R \otimes S$.

\begin{definition}
	For any symmetric monoidal $\infty$-category $\mathbfcal{C}^\otimes$, we define the \emph{$\infty$-category of commutative bialgebra objects in $\mathbfcal C$} to be
	\[ 
		\CBiAlg(\mathbfcal{C}) 
		= \Alg(\CAlg(\mathbfcal{C})^{\op})^{\op}. 
	\]
\end{definition}

	As discussed in \cref{underlying_multiplication}, such a bialgebra produces an object $X \in h\mathbfcal{C}$ with a multiplication map $X \otimes X \to X$, and a refinement to $X \in h\!\CAlg(\mathbfcal{C})$ with a comultiplication map $X \to X \otimes X$, which is well-defined in the homotopy category $h\!\CAlg(\mathbfcal{C})$ and therefore also in the underlying homotopy category $h\mathbfcal C$.
Therefore, for any bialgebra $X$ we can define the \emph{(right) shear map}
\[ 
	\sh\colon X \otimes X \to X \otimes X 
\]
by applying the formula in \eqref{eq:shear_definition} to the multiplication and comultiplication maps in the homotopy category $h\mathbfcal C$.

We say that a bialgebra is a \emph{Hopf algebra} if its shear map is an isomorphism in the homotopy category. 
By \cref{shear_equivalences}, we see that this definition agrees with all the other ones in use.

\begin{definition}\label{df:hopf}
	For any symmetric monoidal $\infty$-category $\mathbfcal{C}^\otimes$, the \emph{$\infty$-category of commutative Hopf algebras in $\mathbfcal{C}$} is the full $\infty$-subcategory
	\[ 
		\CHopf(\mathbfcal{C}) 
			\subseteq 
		\Alg(\CAlg(\mathbfcal{C})^{\op})^{\op} 
	\]
spanned by those objects $X$ for which the shear map \cref{eq:shear_definition} is an isomorphism in $h\mathbfcal C$.
\end{definition}

This notion of Hopf algebra in a symmetric monoidal $\infty$-category is also used by Ergus \cite[Definition 4.0.1]{HopfInfty}. 

\subsection{Translating bialgebras from simplicial categories to \texorpdfstring{$\infty$}{∞}-categories}

Next we record the following lemma, for which we were unable to find a proper reference.

\begin{lemma}\label{lem-sym-mon-pres}
	{Strong} symmetric monoidal functors between symmetric monoidal $\infty$-categories preserve commutative bialgebra and commutative Hopf algebra objects.
\end{lemma}

\begin{proof}
	Let $\mathbfcal{C}$ and $\mathbfcal{D}$ be two symmetric monoidal $\infty$-categories and let $F \colon \mathbfcal{C}^\otimes \rightarrow \mathbfcal{D}^\otimes$ be a {strong} symmetric monoidal functor. 
	We use $F$ to denote both the functor at the level of symmetric monoidal $\infty$-categories and the functor on the underlying $\infty$-categories $\mathbfcal{C} \to \mathbfcal{D}$. 
	We want to show that $F$ lifts to a functor $F_{\CBiAlg}$ on commutative bialgebras that fits into a commuting diagram 
	\begin{center}
		\begin{tikzcd}[column sep = 4em]
			\CBiAlg(\mathbfcal{C}) 
				\arrow[r,"F_{\CBiAlg}"] 
				\arrow[d, "U"] 
				& 
			\CBiAlg(\mathbfcal{D}) 
				\arrow[d, "U"] 
				\\
			\mathbfcal{C} 
				\arrow[r, "F"] 
				& 
			\mathbfcal{D},
		\end{tikzcd}   
	\end{center}
	where the vertical arrows are the forgetful functors.

	As in \cref{induced_functor_of_o_algebras}, we lift $F$ to a functor $F_{\CAlg}$ of commutative algebra objects, simply by post-composing the lax symmetric monoidal functors $N(\Fin_*) \to \mathbfcal{C}^{\otimes}$ with $F$. 
	Since $F$ is {strong} symmetric monoidal, the functor $F_{\CAlg}$ preserves finite coproducts by \cref{sm_coproduct_preserving}.
	Therefore the corresponding functor of opposite categories
	\[ 
		F_{\CAlg}^{\op}\colon \CAlg(\mathbfcal{C})^{\op} 
			\to 
		\CAlg(\mathbfcal{D})^{\op} 
	\]
	is product-preserving, and therefore extends canonically to a {strong} symmetric monoidal functor
	\[ 
		F_{\CAlg}^{\op}\colon (\CAlg(\mathbfcal{C})^{\op})^\times 
			\to 
		(\CAlg(\mathbfcal{D})^{\op})^\times
	\]
	by \cite[2.4.1.8]{HA}. 
	Then by the same reasoning as before, the functor 
	\[
		\begin{tikzcd}[column sep = 6em] 
			(\Fun^{\lax}(\Assoc^\otimes,(\CAlg(\mathbfcal{C})^{\op})^{\times}))^{\op} 
				\rar{(F_{\CAlg}^{\op} \circ -)^{\op}} 
				& 
			(\Fun^{\lax}(\Assoc^\otimes,(\CAlg(\mathbfcal{D})^{\op})^{\times}))^{\op} 
		\end{tikzcd}
	\]
	is also a lift of $F$. 
	Again, by definition this is a functor $\CBiAlg(\mathbfcal{C}) \to \CBiAlg(\mathbfcal{D})$ that we denote $F_{\CBiAlg}$, and we are done with the first claim.

	For the second claim, we recall that {strong} symmetric monoidal functors of $\infty$-categories induce {strong} symmetric monoidal functors on the homotopy category because $F$ preserves cocartesian arrows. 
	Therefore, there are canonical isomorphisms in the homotopy category $F(X) \otimes F(Y) \cong F(X \otimes Y)$ and $\mathbbm{1} \cong F(\mathbbm{1})$. 
	Along these isomorphisms, for a bialgebra $X$, the bialgebra structure on $F(X)$ is obtained by applying $F$ to the multiplication and comultiplication maps of $X$. 
	Therefore the shear map for $F(X)$ is isomorphic in $h\mathbfcal{D}$ to the image of the shear map for $X$. 
	It follows that $F$ preserves commutative Hopf algebras as well.
\end{proof}

Now we are ready to prove \cref{main_comparison-intro}, the main comparison theorem for bialgebras and Hopf algebras.
\begin{itemize}
    \item 
    	Assume $M_\bullet$ is a symmetric monoidal simplicial model category. 
	    We assume that it comes from a topological symmetric monoidal model category, so that the mapping spaces in $M_\bullet$ and $\CAlg(M_\bullet)$ are all Kan complexes (see \cref{top_case_1}).
    \item 
    	Let $M^\flat_\bullet$ be a full subcategory containing the cofibrant objects, closed under the tensor product, and on which the tensor preserves all equivalences. 
	    As before, we let $W^\flat = W \cap M^\flat_0$ be the restriction of the weak equivalences to this subcategory.
\end{itemize}

\begin{theorem}\label{main_comparison}
    Under these assumptions, there is a canonical map of $\infty$-categories
    \[ 
    	N^s(\BiAlg_{(\Comm,D_1)}(M^\flat_\bullet)) 
			\to 
		\CBiAlg(N^s(M^\flat_\bullet)[((W^\flat)^\otimes)^{-1}]) 
			\simeq 
		\CBiAlg(N^s(M^{\mathrm{cf}}_\bullet)) 
	\]
    which preserves the product, coproduct, and shear map as maps in the homotopy category $hM^\flat_\bullet$. 
    In particular, it preserves commutative Hopf algebra objects.
\end{theorem}

The equivalence comes from the equivalence of symmetric monoidal $\infty$-categories from \cref{underlying_sm_inf_cat_compare}, and shows that the right-hand side is the $\infty$-category of commutative bialgebras in the underlying symmetric monoidal $\infty$-category of $M_\bullet$.

\begin{proof}
    First note that there is a canonical isomorphism $N^s(C_\bullet)^{\op} \cong N^s(C_\bullet^{\op})$ for any simplicial category $C_\bullet$.
    
    Next let $D_\bullet = \CAlg(M^\flat_\bullet)^{\op}$, whose mapping spaces are all Kan complexes. 
    The canonical maps we have defined in this paper string together to give a map of $\infty$-categories
    \[ 
    	\begin{array}{rclcl}
            N^s(\BiAlg_{(\Comm,D_1)}(M^\flat_\bullet))
            &= & N^s(\Alg_{D_1}(D_\bullet)^{\op}) & \hspace{1em} & \textup{by definition} 
            	\\[0.3em]
            &\cong & N^s(\Alg_{D_1}(D_\bullet))^{\op} && \textup{$N^s$ commutes with $\op$} 
            	\\[0.3em]
            &\to & \Alg_{\mathbb E_1}(N^s(D_\bullet))^{\op} && \textup{\cref{alg_interchange}} 
            	\\[0.3em]
            &\cong & \Alg(N^s(D_\bullet))^{\op} && \textup{\cref{e1_to_assoc}}.
    	\end{array} 
	\]
    Expanding out the definition of $D_\bullet$ as $\CAlg(M^\flat_\bullet)^{\op}$, we get another string of maps of $\infty$-categories
    \[ 
    	\begin{array}{rclcl}
            N^s(D_\bullet)
            &= & N^s(\CAlg(M^\flat_\bullet)^{\op}) & \hspace{1em} & \textup{by definition} 
            	\\[0.3em]
            &\cong & N^s(\CAlg(M^\flat_\bullet))^{\op} && \textup{$N^s$ commutes with $\op$} 
            	\\[0.3em]
            &\to & \CAlg(N^s(M^\flat_\bullet))^{\op} && \textup{\cref{alg_interchange}} 
            	\\[0.3em]
            &\to & \CAlg(N^s(M^\flat_\bullet)[((W^\flat)^\otimes)^{-1}])^{\op} && \textup{\cref{underlying_sm_inf_cat_compare} and \cref{induced_functor_of_o_algebras}}.
    	\end{array} 
	\]
    Note that by \cref{interchange_coproucts} and a combination of \cref{underlying_sm_inf_cat_compare} and \cref{sm_coproduct_preserving}, the functors on this list are all product-preserving, and are therefore symmetric monoidal in a canonical way by \cite[2.4.1.8]{HA}. 
    Therefore by \cref{induced_functor_of_o_algebras} they induce a map on the $\infty$-category of algebras
    \[ 
    	\Alg(N^s(D_\bullet)) 
			\to 
		\Alg(\CAlg(N^s(M^\flat_\bullet)[((W^\flat)^\otimes)^{-1}])^{\op}). 
	\]
    Taking opposite categories one more time and composing with the first list of functors gives the desired map of $\infty$-categories.
    
    It remains to check that the product, coproduct, and shear map are preserved along the map of homotopy categories $hM^\flat_\bullet \to hM^\flat_\bullet[(W^\flat)^{-1}]$. 
    For the two steps that use the canonical interchange, this follows from \cref{underlying_multiplication}. For the last step where we localize this follows from \cref{lem-sym-mon-pres}, and the fact that when $L$ is the localization functor, the isomorphism $L(X) \otimes L(Y) \cong L(X \otimes Y)$ is just the identity map of $X \otimes Y$.
\end{proof}

 \section{An application}

As an application, we illustrate one method for applying \cref{main_comparison} to produce a Hopf algebra spectrum in the $\infty$-category of spectra. We use this method in \cite{KKMMW1} to endow reduced spherical scissors congruence $K$-theory with the structure of an $(E_{\infty},E_1)$-Hopf algebra spectrum.

Recall that a theorem of P\'{e}roux--Shipley precludes the existence of non-cocommutative coalgebras in any symmetric monoidal model category of spectra \cite{PS_coalgebras}. To circumvent this obstacle, we construct a Hopf algebra in a $1$-category of diagrams of spectra, and then pass to the underlying $\infty$-category of diagrams of spectra. Taking the colimit of the diagram yields a Hopf algebra in the $\infty$-category of spectra.

Let $\Sp^O_\bullet$ be the symmetric monoidal topological model category of orthogonal spectra with the stable model structure \cite[Theorem 9.2]{mandell_may_shipley_schwede}. We consider this as a symmetric monoidal simplicial model category via the singular simplices functor. We write $\mathbf{Sp}$ for the $\infty$-category of spectra, which is the underlying symmetric monoidal $\infty$-category of $\Sp^O_\bullet$.

Let $I$ be any symmetric monoidal $1$-category. The $1$-category of diagrams $\Fun(I,\Sp^O_0)$ then becomes a simplicial model category $\Fun(I,\Sp^O_\bullet)$ with the projective model structure in which weak equivalences and fibrations are determined pointwise \cite[A.3.3.1]{HTT}. As in \cref{too many bullets}, we abuse notation and suppress the additional placeholder index that would be necessary on $\Fun(I, \Sp^O_\bullet)$ to indicate that it is a simplicial category. Furthermore, this simplicial model category is a symmetric monoidal simplicial model category via Day convolution \cite[Theorem 5.6.35]{purplebook} \cite[Proposition 8.4(2)]{isaacson}, using the symmetric monoidal structure of $I$ and the smash product on $\Sp^O_\bullet$. We let $\Fun(I,\Sp^O_\bullet)^\flat$ be the full subcategory of the pointwise cofibrant diagrams, ie those diagrams $F \colon I \to \Sp^O_\bullet$ such that $F(i)$ is cofibrant for each $i \in I$. This strictly contains the cofibrant diagrams in the projective model structure.

\begin{corollary}\label{main_corollary}
	Under these assumptions, there is a canonical map of $\infty$-categories
	\[ 
		N^s(\BiAlg_{(\Comm,D_1)}(\Fun(I,\Sp^O_\bullet)^\flat)) \to  \CBiAlg(\mathbf{Sp})
	\]
    that sends each diagram $X$ to its homotopy colimit. Furthermore this preserves the product, coproduct, and shear map as maps in the homotopy category, and therefore also preserves Hopf algebra objects.
\end{corollary}

In particular, for every diagram $X\colon I \to \Sp^O_0$ that is pointwise cofibrant and a $(\Comm,D_1)$-bialgebra in the 1-category of diagrams of spectra, the homotopy colimit of $X$ naturally has the structure of a commutative bialgebra in the $\infty$-category of spectra.

\begin{proof}
\cref{main_comparison} gives us a map of $\infty$-categories
\begin{equation}
\label{main comparison in functor category}
	N^s(\BiAlg_{(\Comm,D_1)}(\Fun(I,\Sp^O_\bullet)^\flat)) \to \CBiAlg(N^s(\Fun(I, \Sp^O_\bullet)^{cf}))
\end{equation}
which preserves Hopf algebra objects. All we need to do is compose this with the homotopy colimit functor, but for this to work, we need to verify that the homotopy colimit induces a functor on the $\infty$-categories of commutative bialgebras.

We define the homotopy colimit to be the left-derived functor of the colimit functor of 1-categories
\[
	\colim_I \colon \Fun(I, \Sp^O_0) \to \Sp^O_0.
\]
This is left Quillen with respect to the projective model structure on $\Fun(I, \Sp^O_0)$, and is strong symmetric monoidal as a functor between $1$-categories as in e.g. \cite{ben-moshe-schlank}. Therefore if we restrict to cofibrant diagrams in the projective model structure and take the nerve, we get a strong symmetric monoidal functor of $\infty$-categories
\[
	N(\Fun(I, \Sp^O_0)^c) \to N((\Sp^O_0)^c).
\]
The weak equivalences are all preserved by this functor and by the tensor products, so by the universal property from \cref{sm_localization}, inverting the weak equivalences gives a strong symmetric monoidal functor of $\infty$-categories
\[
	N(\Fun(I, \Sp^O_0)^c)[(W^c)^{-1}] \to N((\Sp^O_0)^c)[(W^c)^{-1}].
\]
Using the standard equivalence $N(M^c_0)[(W^c)^{-1}] \simeq N(M_0)[W^{-1}]$ and the equivalences in \cref{underlying_sm_inf_cat_compare_2}, we can rewrite this as a strong symmetric monoidal functor
\[ N^s(\Fun(I, \Sp^O_\bullet)^{cf}) \to N^s((\Sp^O_\bullet)^{cf}) \simeq \mathbf{Sp}. \]
By \cref{lem-sym-mon-pres}, this induces a functor on the $\infty$-categories of bialgebras
\begin{equation}
\label{main comparison in functor category 2}
    \CBiAlg(N^s(\Fun(I, \Sp^O_\bullet)^{cf})) \to \CBiAlg(N^s((\Sp^O_\bullet)^{cf})) \simeq \CBiAlg(\mathbf{Sp})
\end{equation}
which preserves Hopf algebra objects. Now we compose \eqref{main comparison in functor category} with \eqref{main comparison in functor category 2} to get the desired map of $\infty$-categories.
\end{proof}

In \cite{KKMMW1}, we define a version of reduced spherical scissors congruence $K$-theory, which forms a pointwise cofibrant commutative Hopf algebra $\widetilde{\mathfrak{S}}$ in the symmetric monoidal simplicial model category $\Fun(I, \Sp^O_\bullet)$: 
\[
	\widetilde{\mathfrak{S}} \in \BiAlg_{(\Comm,D_1)}(\Fun(I,\Sp^O_\bullet)^\flat).
\]
By \cref{main_corollary}, its homotopy colimit $\widetilde{\mathcal{S}}$ is therefore a commutative Hopf algebra in the $\infty$-category of spectra.

\bibliographystyle{alpha}
\bibliography{references.bib}
\end{document}